\documentclass[times]{nlaauth}
\usepackage{enumerate}
\usepackage{moreverb}
\usepackage{amsthm}
\usepackage{amsmath}
\usepackage{graphicx}
\usepackage{tabularx}
\usepackage{subfigure}
\usepackage{mathrsfs}
\usepackage{multirow}
\usepackage{pgfplots}

\newcommand\BibTeX{{\rmfamily B\kern-.05em \textsc{i\kern-.025em b}\kern-.08em
T\kern-.1667em\lower.7ex\hbox{E}\kern-.125emX}}
\newcommand{\diag}{\mathop{\mathrm{diag}}\nolimits}

\newtheorem{theorem}{Theorem}

\newcounter{algorithm}
\newenvironment{algorithm}{\refstepcounter{algorithm}\vspace{1ex}
  {\sc Algorithm \thealgorithm.}\hspace{0.3em}\parindent=0pt}{\vspace{1ex}}
\newcounter{remark}
\newenvironment{remark}{\refstepcounter{remark}\vspace{1ex}
{\sc Remark \theremark.}\hspace{0.3em}\parindent=0pt}{\vspace{1ex}}

\def\volumeyear{2015}

\begin{document}

\runningheads{Tridiagonal eigenvalue problem}{HSS matrix}

\title{New fast divide-and-conquer algorithms for the symmetric tridiagonal eigenvalue problem}

\author{Shengguo~Li\corrauth, Xiangke Liao, Jie Liu, and Hao Jiang}

\address{College of Computer Science, National University of Defense
  Technology, Changsha 410073, China}

\corraddr{nudtlsg@nudt.edu.cn}

\begin{abstract}
  In this paper, two accelerated divide-and-conquer algorithms are proposed for
  the symmetric tridiagonal eigenvalue problem, which cost $O(N^2r)$ {flops}
  in the worst case, where $N$ is the dimension of the matrix and $r$ is a modest number depending
  on the distribution of eigenvalues.
  Both of these algorithms use hierarchically semiseparable (HSS) matrices to approximate some
  intermediate eigenvector matrices which are Cauchy-like matrices and are off-diagonally low-rank.
  The difference of these two versions lies in using different HSS construction algorithms,
  one (denoted by {ADC1}) uses a structured low-rank approximation method
  and the other ({ADC2}) uses a randomized HSS construction algorithm.
  For the ADC2 algorithm, a method is proposed to estimate the off-diagonal rank.
  Numerous experiments have been done to show their stability and efficiency.
  These algorithms are implemented in parallel in a shared memory environment, and
  some parallel implementation details are included. Comparing the ADCs with
  highly optimized multithreaded libraries such as Intel MKL, we find that ADCs could be more than 6x times faster
  for some large matrices with few deflations.
\end{abstract}


\maketitle


\vspace{-6pt}

\section{Introduction}
\vspace{-2pt}

Computing the eigendecomposition of a symmetric tridiagonal matrix is a classic linear algebra
problem, and is ubiquitous in computational science.
Some well-known algorithms include the QR algorithm~\cite{Golub-book2,Parlett-book},
the MRRR~\cite{MRRR-LAA} and the divide-and-conquer (DC) algorithm~\cite{Cuppen81,Gu-eigenvalue}.
According to the comparisons in~\cite{Demmel-eig-perf}, DC and MRRR are generally faster than QR for large matrices, especially when
the eigenvectors are required.
In this work, we focus on the DC algorithm, and the goal is to develop an improved version.

Though DC is very fast in practice which takes $O(N^{2.3})$ flops on average~\cite{Demmel-book,parallel_dc},
for some matrices with few deflations its complexity \cite{Cuppen81,Rutter94,Gu-eigenvalue} can be $O(N^3)$.
By using hierarchically semiseparable (HSS) matrices~\cite{Hss-ulv,Xia-Fast09},
we show that its worst case complexity can be reduced to $O(N^2r)$, where $r$ is a modest
number and is usually much smaller than a big $N$.
The main observation of this work is from the following theorem.

\begin{theorem}[Bunch, Nielsen, and Sorensen~\cite{BNS-Rankone}]
  \label{thm:BNS-R1}
 Assume that $D=\diag(d_1,\cdots,d_N)$ such that $d_1<d_2<\cdots<d_N$, and $u\in \mathbb{R}^N$ is a vector and a scalar $\rho > 0$. Then, the eigenvector corresponding to $\lambda_i$, an eigenvalue of $M=D+\rho uu^T$, is
  \begin{equation}
    \label{eq:EVec}
    q_i =
    \begin{pmatrix}
      \frac{u_1}{d_1-\lambda_i},\cdots,\frac{u_N}{d_N-\lambda_i}
    \end{pmatrix}^T /
    \sqrt{\sum_{j=1}^N \frac{u_j^2}{(d_j-\lambda_i)^2}},
  \end{equation}
  and the eigenvalues of $M$ satisfy
  \begin{equation}
    \label{eq:ev-interlacing}
    d_1 < \lambda_1 < d_2 < \lambda_2 < \cdots < d_N < \lambda_N.  \hspace{1.3in} \square
  \end{equation}
\end{theorem}

It shows that $Q=\left( \frac{u_i v_j}{d_i-\lambda_j} \right)_{i,j}$ with
$v_j=1/\sqrt{\sum_{k=1}^N \frac{u_k^2}{(d_k-\lambda_j)^2}}$ is a Cauchy-like matrix.
Recall that $C$ is called Cauchy-like if it satisfies
\begin{equation}
\label{eq:cauchylike}
  D\cdotp C-C\cdotp \Lambda=u\cdotp v^T,
\end{equation}
  where $D=\diag(d_1,\cdots, d_N), \Lambda=\diag(\lambda_1,\cdots,\lambda_N)$ and $u, v\in \mathbb{R}^N$,
  which are called the \emph{generators} of Cauchy-like matrix $Q$.
  It is easy to check that $Q$ is also off-diagonally low-rank.
  To take advantage of these two properties, we can use an HSS matrix to approximate $Q$ and then
  use the fast HSS matrix multiplication algorithm to update the eigenvectors, like the bidiagonal SVD case~\cite{Shengguo-SIMAX2}.
  A structured low-rank approximation method is designed for
  a Cauchy-like matrix in~\cite{GX-Toeplitz,Shengguo-SIMAX2}, called SRRSC (\emph{structured rank-revealing
  Schur-complement factorization}), which can be used to construct HSS matrices efficiently.
  By incorporating SRRSC into DC, an accelerated DC (ADC) algorithm is proposed for the singular value problem in~\cite{Shengguo-SIMAX2},
  where ADC is 3x faster than DC in Intel MKL for some large matrices.
  We show that this technique also works for the symmetric tridiagonal eigenvalue problem,
  which is presented in Example 3 in section~\ref{sec:numer-tedc}.
  In this paper, we implement the ADC algorithm in parallel and it obtains even better speedups.

  In this paper, we further show that the randomized HSS construction algorithm~\cite{rand-hss} can
  also be used to compute the eigendecomposition reliably, and it can
  obtain similar speedups as using SRRSC when compared with Intel MKL.
  Recent work has suggested the efficiency of the randomized algorithms~\cite{Martinsson-Rev10,Martinsson-PNAS07} for computing a low-rank
  matrix approximation.
  Martinsson~\cite{rand-hss} has developed a novel HSS construction algorithm by using the randomized sampling technique, and
  for simplicity we refer to this algorithm as \texttt{RSHSS}.
  This method is extremely suitable for matrices with fast matrix-vector multiplication algorithms.
  For example, if the matrix-vector product costs $O(N)$ flops, RSHSS
  has linear complexity $O(Nr^2)$, where $N$ is the dimension of the matrix and $r$ is its maximum numerical rank of off-diagonal blocks.
  For matrix $Q$ in~\eqref{eq:EVec}, the fast multipole method (FMM)~\cite{FMM87,Rokhlin88} can be used for the
  matrix-vector product in $O(N \log N)$ flops.
  Therefore, an HSS matrix approximation to $Q$ can be constructed in $O(Nr^2)$ flops if combining RSHSS with FMM.
  To use RSHSS, we need to know the maximum rank of off-diagonal blocks,
  which is difficult for general matrices. Fortunately,
  the off-diagonal rank of the matrix $Q$ defined in~\eqref{eq:EVec} can be estimated by using the approximation theory of function $1/x$.
  We show in section~\ref{sec:notation} that the estimated rank based on the exponential expansion~\cite{hackbusch-exp}
  is quite acceptable.

  By adding the HSS matrix techniques to the symmetric DC algorithm~\cite{Cuppen81,Rutter94,Gu-eigenvalue},
  two new accelerated DC (ADC) algorithms are obtained. One is denoted by {ADC1} based on SRRSC,
  and the other is denoted by {ADC2} based on RSHSS.
  Similar to the analysis in~\cite{Demmel-book}, the complexity of
  ADCs can be shown to be $O(N^2 r)$ where $N$ is the dimension of a symmetric tridiagonal matrix $T$
  and $r$ is a modest integer, which is related to the distribution of eigenvalues.
  Since the HSS matrix construction~\cite{Hss-ulv,rand-hss,Sherry-PHss}
  and multiplication~\cite{Lyons-thesis} algorithms are naturally parallelizable,
  we further implement these two algorithms in parallel by using OpenMP in a shared
  memory multicore environment.
  We also simply parallelize the classical process of DC algorithm by using OpenMP such
  as solving the subproblems at the bottom level of the divide-and-conquer tree and all the secular equations.
  Numerous experiments have been done to test these two ADC algorithms.
  It turns out that our ADCs can be about 6x times faster than the DC implementation in MKL for some large matrices with few deflations.
  The accuracy comparisons are also included in section~\ref{sec:numer-tedc}.

\section{Preliminary}
\label{sec:notation}

Assume that $T$ is a symmetric tridiagonal matrix,
\begin{equation}
\label{eq:T}
T=\text{tridiag}\left( \begin{array}{ccccccccccc}
& b_1 & & b_2 & & \cdotp & & b_{N-2} & & b_{N-1} & \\
a_1 & & a_2 & & \cdotp & & \cdotp & &  a_{N-1} & & a_N \\
& b_1 & & b_2 & & \cdotp & & b_{N-2} & & b_{N-1} &
 \end{array} \right).
\end{equation}
We briefly introduce some formulae of Cuppen's divide-and-conquer
algorithm~\cite{Cuppen81,Rutter94}.
$T$ is decomposed into the sum of two matrices,
\begin{equation}
  \label{eq:T2}
  T=
  \begin{bmatrix}
    T_1 & \\ & T_2
  \end{bmatrix}+b_k vv^T,
\end{equation}
where $T_1\in R^{k\times k}$ and $v=[0,\ldots,0,1,1,0,\ldots,0]^T$ with
ones at the $k$-th and $(k+1)$-th entries.
If $T_1=Q_1 D_1Q_1^T$ and $T_2=Q_2 D_2 Q_2^T$, then $T$ can be written as
\begin{equation}
  \label{eq:T3}
  T=
  \begin{bmatrix}
    Q_1 & \\ & Q_2
  \end{bmatrix} \left(
  \begin{bmatrix}
    D_1 & \\ & D_2
  \end{bmatrix} + b_k uu^T \right)
\begin{bmatrix}
  Q_1^T & \\ & Q_2^T
\end{bmatrix},
\end{equation}
where $u=
\begin{bmatrix}
  Q_1^T & \\ & Q_2^T
\end{bmatrix} v =
\begin{bmatrix}
  \text{last column of } Q_1^T \\
  \text{first column of } Q_2^T
\end{bmatrix}.
$
By Theorem~\ref{thm:BNS-R1}, we can get the eigenvectors $Q$ of the middle
matrix at the right hand side of~\eqref{eq:T3}.
The eigenvectors of $T$ would be $\begin{bmatrix} Q_1 & \\ & Q_2\end{bmatrix}Q$.

The matrices $T_i, i=1,2$ can also be divided recursively.
More numerical details can be found in~\cite{Rutter94} and \cite{Demmel-book}. We only point
out that the eigenvectors can not be computed directly from~\eqref{eq:EVec}.
In practice, the computed $\hat{\lambda}_i$ is only an approximation to $\lambda_i$.
To compute the eigenvectors orthogonally~\cite{Gu-eigenvalue,Gu-rank1},
we need to use L\"{o}wner's Theorem~\cite{lowner} to recompute the vector $u$ as
\begin{equation}
  \label{eq:ev_z}
  \hat{u}_i = \sqrt{\prod_{j=1}^{i-1}\left(\frac{\hat{\lambda}_j-d_i}{d_j-d_i}\right)\cdotp
    \prod_{j=i+1}^{N}\left(\frac{\hat{\lambda}_j-d_i}{d_{j}-d_i}\right) \cdotp (\hat{\lambda}_i-d_i)},
\end{equation}
and then use~\eqref{eq:EVec} to compute the eigenvectors.

\subsection{The low-rank structure of $Q$}
\label{sec:lowrank}

 To be more specific, the matrix $Q$ is defined as
 \begin{equation}
   \label{eq:Ev-Q}
   Q=
   \begin{bmatrix}
     \frac{u_1 v_1}{d_1 - \lambda_1} & \frac{u_1 v_2}{d_1-\lambda_2} & \cdots & \frac{u_1 v_N}{d_1-\lambda_N} \\
     \frac{u_2 v_1}{d_2 - \lambda_1} & \frac{u_2 v_2}{d_2-\lambda_2} & \cdots & \frac{u_2v_N}{d_2-\lambda_N} \\
     \vdots & \vdots & & \vdots \\
     \frac{u_N v_1}{d_N - \lambda_1} & \frac{u_N v_2}{d_N-\lambda_2} & \cdots & \frac{u_N v_N}{d_N-\lambda_N} \\
   \end{bmatrix}.
 \end{equation}
Since $\{d_i\}$ and $\{\lambda_j\}$ are interlacing, see~\eqref{eq:ev-interlacing}, $Q$ is usually
off-diagonally low rank, and the ranks of off-diagonal blocks depend on the distribution of $\{d_i\}$ and $\{\lambda_i\}$.
We use the following example to show that.

{\bf Example 1.}
Assume that a matrix $Q$ satisfies the eigendecomposition $M=D+uu^T=Q\Lambda Q^T$,
where $D=\diag(d_1,d_2,\ldots,d_N)$, $d_i=i\cdotp \frac{b-a}{N}$, $a=1.0$, $b=9.0$,
for $i=1,\ldots,N$ and
$u\in \mathbb{R}^N$ is a random normalized vector.
The off-diagonal low rank property of $Q$ is shown in Table~\ref{tab:Ex1-rank},
which includes the numerical ranks of the submatrices $Q(1:m,m+1:N)$ for different $m$ with $N=2000$.
The ranks are computed by truncating the singular values less than $1.0e^{-13}$.

\begin{table}[ptbh]
\caption{The ranks of different off-diagonal blocks of $Q$}
\label{tab:Ex1-rank}
\begin{center}
\begin{tabular}[c]{|c|cccccccccc|} \hline
$m$ & 100 & 200 & 300 & 400 & 500 & 600 & 700 & 800 & 900 & 1000 \\ \hline
\emph{rank} & 18 & 20 & 21 & 22 & 23 & 23 & 23 & 24 & 24 & 24  \\ \hline
\end{tabular}
\end{center}
\end{table}

The ranks of the off-diagonal blocks can be estimated by using the
approximation theory of function $f(x)=1/x$.
The element $\frac{1}{d_i-\lambda_j}$ can be approximated by the sums of exponentials,
\begin{equation}
\label{eq:exponential}
\frac{1}{d_i - \lambda_j} \approx \sum_{k=1}^r \omega_k e^{-\alpha_k (d_i-\lambda_j)}:=s_r(d_i-\lambda_j).
\end{equation}
Assume that $d_i$ and $\lambda_j$ belong to two different subintervals of $[\hat{a}, \hat{b}]$, $d_i \in I_1$, $\lambda_j \in I_2$, $I_1 \cap I_2=\emptyset$ and $d_i <\lambda _j$.
Denote $\mbox{dist}(I_1,I_2)=\min_{d_i\in I_1, \\ \lambda_j \in I_2} |d_i-\lambda_j|$, (In our case $d_i$ and $\lambda_j$
are the eigenvalues of $D$ and $M$ respectively, $\hat{a}=d_1$ and $\hat{b}=\lambda_N$, the largest eigenvalue of $M$), then
\begin{equation}
\label{eq:dist_xy}
1\leq \frac{d_i-\lambda_j}{\mbox{dist}(I_1,I_2)} \le \frac{\hat{b}-\hat{a}}{\mbox{dist}(I_1,I_2)}:=R.
\end{equation}
An approximation error bound is given in~\cite{hackbusch-exp} for the sums of exponentials,
\begin{equation}
\label{eq:app}
\left| \frac{1}{x}-s_r(x) \right| \leq 16 e^{-\frac{r\pi^2}{\log(8R)}},
\end{equation}
where $s_r(x)$ is defined in~\eqref{eq:exponential} and $x\in [1, R]$.
As long as the number of approximation terms satisfies
\begin{equation}
\label{eq:bounds}
r \ge \lceil \frac{\log(16/\epsilon) \log(8R)}{\pi^2} \rceil,
\end{equation}
the approximation error of the sum of exponentials is less than $\epsilon$, a small constant.
For example, when $m=300$ and $\epsilon=1e$-$13$, $\mbox{dist}(I_1,I_2)=4.3e$-$3$, $R=2.0e3$ and
the off-diagonal rank estimated by~\eqref{eq:app} is 32, which is close to the result in
Table~\ref{tab:Ex1-rank}.
Equation~\eqref{eq:app} would be used to estimate the rank in Algorithm~\ref{alg:randhss}.

\begin{remark}
To keep the orthogonality of $Q$ , $d_i-\lambda_j$ must be computed by
\begin{equation}
  \label{eq:d-lambda}
  d_i - \lambda_j =
\begin{cases}
(d_i-d_j)-\gamma_j & \text{ if } i\le j  \\
(d_i -d_{j+1})+\mu_j & \text{ if } i > j
\end{cases},
\end{equation}
where $\gamma_i=\lambda_i-d_i$ (the distance between $\lambda_i$ and $d_i$),
and $\mu_i=d_{i+1}-\lambda_i$ (the distance between $\lambda_i$ and $d_{i+1}$),
which can be returned by calling the LAPACK routine \texttt{dlaed4}.
If using FMM to compute the eigenvectors, equation~\eqref{eq:d-lambda} shows that some modifications of classic FMM are needed
since some $d_i$ and $\lambda_j$ may equal in double precision but
$\gamma_i$ is not zero.

\end{remark}

\subsection{Introduction to HSS matrices}
\label{sec:hss}

The HSS matrices are a very important type of rank-structured matrices
which share the same property that the off-diagonal blocks are low-rank.
Other rank-structured matrices include  $\mathcal{H}$-matrix \cite{Hackbusch1999,Hackbusch2000},
$\mathcal{H}^{2}$-matrix \cite{Hackbusch-Sauter2000,Hackbusch-Borm2002},
quasiseparable matrices \cite{Eidelman-Gohberg1999, Vandebril-book1}, and sequentially
semiseparable (SSS) \cite{Chandrasekaran03,ChandrasekaranGu05} matrices.
The HSS matrix was first discussed in~\cite{ChandrasekaranGu04,Hss-ulv}, which can be seen as
an algebraic counterpart of FMM in~\cite{Starr-thesis}.
In this paper we use the HSS matrix to accelerate the computation of eigenvectors,
and other rank-structured matrices can be similarly used too.
We follow the notation in~\cite{Xia-Fast09,Xia-HSS-Chol} and briefly introduce some key concepts of
the HSS matrix.

Let $\mathcal{I}=\{1,2,\ldots, N\}$ and
$\mathcal{T}$ be a \emph{postordered} binary tree, which means the ordering of a nonleaf node $i$
satisfies $i_1 < i_2 <i$, where $i_1$ is its left child and $i_2$ is its right child.
Each node $i$ is associated with a contiguous subset of $\mathcal{I}$, $t_i$, satisfying the
following conditions:
\begin{itemize}
\item $t_{i_1}\cup t_{i_2}=t_i$ and $t_{i_1}\cap t_{i_2}=\emptyset$, for a parent node $i$ with
  left child $i_1$ and right child $i_2$;

  \item $\cup_{i \in LN} t_i = \mathcal{I}$, where $LN$ denotes the set of all leaf nodes;

\item $t_{\operatorname*{root}(\mathcal{T})}=\mathcal{I}$,
  $\operatorname*{root}(\mathcal{T})$ denotes the root of $\mathcal{T}$.
\end{itemize}
A block row or column excluding the diagonal block is called an \emph{HSS block row or column},
denoted by
\[
H_i^{row} = A_{t_i \times (\mathcal{I}\backslash t_i)}, \quad H_i^{col} = A_{(\mathcal{I}\backslash t_i)\times t_i},
\]
associated with node $i$. We also simply call them \emph{HSS blocks}.
As in~\cite{Xia-HSS-Chol}, the maximum (numerical) rank of all the HSS blocks is
called \emph{HSS rank}.

\begin{figure}
\centering
\subfigure[Matrix $A$]{
\includegraphics[width=1.8in]{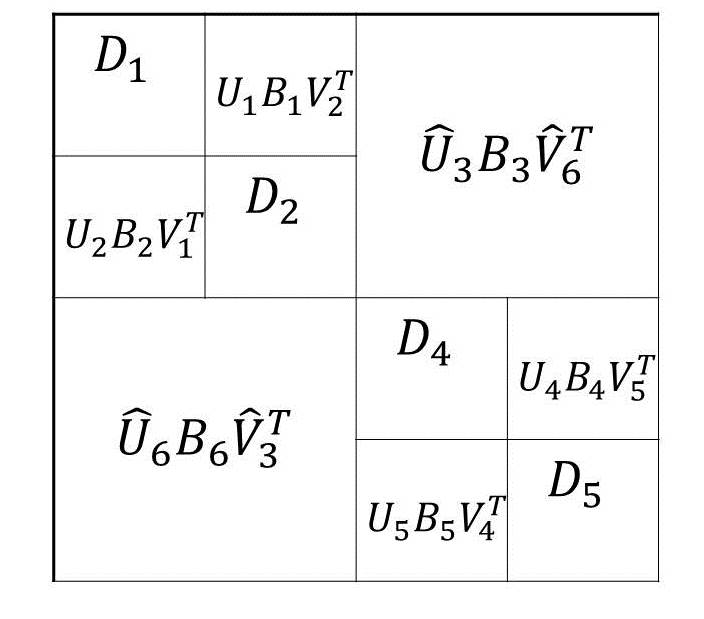}
\label{fig:Ahss}}
\quad
\subfigure[HSS tree]{
\includegraphics[width=1.6in]{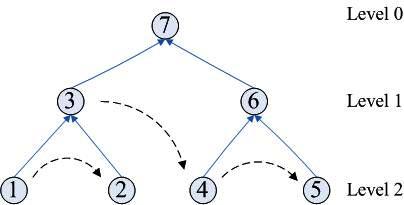}
\label{fig:Hsstree}}
\caption{A $4\times 4$ HSS matrix and its corresponding HSS tree}
\label{fig:gnaHss}
\end{figure}

For each node $i$ in $\mathcal{T}$, there are matrices $\widehat{D}_i$, $\widehat{U}_i$, $\widehat{V}_i$ and $B_i$
associated with it, called \emph{generators}, such that
\begin{equation}
  \label{eq:new-hss}
  \begin{split}
  \widehat{D}_i & =A|_{t_i\times t_i}=
  \begin{bmatrix} \widehat{D}_{i_1} & \widehat{U}_{i_1}B_{i_1}\widehat{V}_{i_2}^T \\
    \widehat{U}_{i_2} B_{i_2} \widehat{V}_{i_1}^T & \widehat{D}_{i_2} \end{bmatrix}, \\
  \widehat{U}_i & =
  \begin{bmatrix}
    \widehat{U}_{i_1} & \\ & \widehat{U}_{i_2}
  \end{bmatrix}
  U_i, \quad
  \widehat{V}_i =
  \begin{bmatrix}
    \widehat{V}_{i_1} & \\ & \widehat{V}_{i_2}
  \end{bmatrix} V_i.
  \end{split}
\end{equation}
For a leaf node $i$, $\widehat{D}_i=D_i$, $\widehat{U}_i = U_i$,
$\widehat{V}_i=V_i$. Figure~\ref{fig:Ahss} shows a $4\times 4$ HSS matrix $A$, and
it can be written as
\begin{equation}
\label{eq:posthsslevel}
A=\begin{bmatrix} \begin{bmatrix} D_1 & {U}_1B_1{V}_2^T \\  {U}_2B_2{V}_1^T & D_2 \end{bmatrix}   & \widehat{U}_3B_3\widehat{V}_6^T \\
\widehat{U}_6B_6\widehat{V}_3^T &
\begin{bmatrix} D_4 & {U}_4B_4{V}_5^T \\  {U}_5B_5{V}_4^T & D_5 \end{bmatrix}  \end{bmatrix},
\end{equation}
and Figure~\ref{fig:Hsstree} shows its corresponding postordering HSS tree.

\begin{remark}
\begin{enumerate}
  \item The generators of a Cauchy-like matrix~\eqref{eq:cauchylike} can be represented by four vectors.
  While, the generators of an HSS matrix are \emph{matrices}.
  For an HSS matrix, we only need to store the generators $D_i$, $U_i$, $V_i$ and $B_i$, and
  $\widehat{D}_i$, $\widehat{U}_i$ and $\widehat{V}_i$ can be constructed hierarchically when needed.

  \item The HSS representation~\eqref{eq:new-hss} is equivalent to the previous representations in~\cite{Hss-ulv,ChandrasekaranGu04,Xia-HSS-Chol},
   but is simpler (generators $R$ and $W$ are not introduced).
   For a parent node $i$, if let $U_i=\begin{bmatrix} R_{i_1} \\ R_{i_2} \end{bmatrix}$, $V_i=\begin{bmatrix} W_{i_1} \\ W_{i_2} \end{bmatrix}$, then
   \eqref{eq:posthsslevel} reduces to the form in~\cite{Xia-HSS-Chol},
   \begin{equation*}
     A=\left[
\begin{array}
[c]{cccc}%
D_{1} & U_{1}B_{1}V_{2}^{T} & U_{1}R_{1}B_{3}W_{4}^{T}V_{4}^{T} & U_{1}%
R_{1}B_{3}W_{5}^{T}V_{5}^{T}\\
U_{2}B_{2}V_{1}^{T} & D_{2} & U_{2}R_{2}B_{3}W_{4}^{T}V_{4}^{T} & U_{2}%
R_{2}B_{3}W_{5}^{T}V_{5}^{T}\\
U_{4}R_{4}B_{6}W_{1}^{T}V_{1}^{T} & U_{4}R_{4}B_{6}W_{2}^{T}V_{2}^{T} & D_{4}
& U_{4}B_{4}V_{5}^{T}\\
U_{5}R_{5}B_{6}W_{1}^{T}V_{1}^{T} & U_{5}R_{5}B_{6}W_{2}^{T}V_{2}^{T} &
U_{5}B_{5}V_{4}^{T} & D_{5}%
\end{array}
\right]. \label{eq:posterHSS}%
\end{equation*}

\end{enumerate}
\end{remark}

\subsection{Accelerated tridiagonal DC algorithm}
\label{sec:sdc}
The procedure of ADC algorithms is expressed in the following algorithm,
which is similar to the DC algorithm~\cite{Demmel-book}.

\begin{algorithm}[\texttt{ADC}($T,Q,\Lambda$)]
\label{alg:tedc}
Compute the whole eigendecomposition of a symmetric tridiagonal matrix by using the ADC algorithm. Let $m(=25)$ be a small integer constant.

\begin{description}
  \item \texttt{if} the row dimension of $T$ is less than $m$

  \item \hspace{0.5cm} use QR algorithm~\cite{Parlett-book} to compute $T=Q\Lambda Q^T$;
  \item \hspace{0.5cm} \texttt{return} $Q$ and $\Lambda$;

  \item \texttt{else}
  \item \hspace{0.5cm} form $T=\begin{bmatrix} T_1 & \\ & T_2 \end{bmatrix}+b_kvv^T$;
  \item \hspace{0.5cm} call \texttt{ADC}($T_1,Q_1,\Lambda_1$);
  \item \hspace{0.5cm} call \texttt{ADC}($T_2,Q_2,\Lambda_2$);
  \item \hspace{0.5cm} form $M=D+b_kuu^T$ from $Q_1,Q_2,\Lambda_1,\Lambda_2$;
  \item \hspace{0.5cm} \texttt{if} the size of $M$ is small
  \item \hspace{1.0cm} find the eigenvalues $\Lambda$ and eigenvectors $Q^\prime$ of $M$;
  \item \hspace{1.0cm} compute $Q=\begin{bmatrix} Q_1 & \\ & Q_2 \end{bmatrix}\cdotp Q^\prime$;
  \item \hspace{0.5cm} \texttt{else}
  \item \hspace{1.0cm} find the eigenvalues $\Lambda$ of $M$ and construct an HSS matrix $H_Q\approx Q^\prime$;
  \item \hspace{1.0cm} compute $Q=\begin{bmatrix} Q_1 & \\ & Q_2 \end{bmatrix}\cdotp H_Q$ via the HSS matrix multiplication algorithm;
  \item \hspace{0.5cm} \texttt{end if}
  \item \hspace{0.5cm} \texttt{return} $Q$ and $\Lambda$;
  \item \texttt{end if}
\end{description}

\end{algorithm}

\begin{remark}
Algorithm~\ref{alg:tedc} is more like a framework for accelerating the tridiagonal DC algorithm,
since the HSS matrices can be replaced by other rank-structured matrices such as $\mathcal{H}$-,$\mathcal{H}^2$-matrix.
The difference between Algorithm~\ref{alg:tedc} and the standard DC algorithm is that
Algorithm~\ref{alg:tedc} uses the HSS matrix techniques to update the eigenvectors when
the size of matrix $M$ is large, to reduce the complexity cost.
If the sizes of secular equations are always small, i.e., most eigenvalues are computed by deflation, ADC is equivalent to the standard DC.
\end{remark}

\section{Randomized HSS construction algorithm}

The random HSS construction algorithm proposed in~\cite{rand-hss} is built on two low-rank approximation algorithms:
\emph{random sampling} (RS)~\cite{Martinsson-Rev10,Martinsson-PNAS07}
and \emph{interpolative decomposition} (ID)~\cite{Cheng-Random,Martinsson-Harmon}.
The form of ID has appeared in the rank-revealing QR~\cite{Gu-RRQR} and rank-revealing LU factorization~\cite{Pan00}, and
it also has a close relationship to the matrix skeleton and CUR factorization~\cite{GTZ-ACA,Stewart-QR}.

We introduce RS first. For a given $m\times n$ matrix $B$ with $m < n$, we want to find a tall matrix
$Q$ with orthogonal columns such that
\[
\|B - QQ^* B\| < \epsilon,
\]
where $\epsilon$ is a small constant.
The random sampling method right multiplies $B$ with a Gaussian random matrix $\Omega \in \mathbb{R}^{n\times (r+p)}$,
and get a ``compressed'' matrix $Y=B\Omega$ with much fewer columns, $(r+p) \ll n$,
where $r$ is the numerical rank of $B$ and $p$ is the oversampling parameter, usually $p=5, 10$ or $20$.
Then, the matrix $Q$ can be obtained by applying the RRQR~\cite{Chan92,Gu-RRQR} or the truncated SVD~\cite{Golub-book2} to $Y$.
It is shown in~\cite{rand-hss,Martinsson-Rev10} that the RS algorithm computes a good low-rank
approximation with quite high probability. For example, the computed $Q$ satisfies
\[
\|B-QQ^TB\| \le \left( 1+11\sqrt{(r+p)\min(m,n)} \right)\sigma_{r+1},
\]
with probability at least $1-6p^{-p}$~\cite{rand-hss}.

\begin{remark}
In general, the rank $r$ is rarely known in advance. For the symmetric tridiagonal DC algorithm, we can
use~\eqref{eq:app} as a guide to estimate $r$.
\end{remark}

The ID method computes an approximate low-rank factorization of $B$ such that
\[
B \approx B(:,J) \widetilde{X} \cdotp P =B(:,J)X,
\]
where $J$ is a subset of the column indices of $B$, $\widetilde{X}$ is a $r\times n$ matrix with
a $r\times r$ identity matrix as a submatrix and all its entries are less than one in
magnitude, and $P$ is a permutation matrix.
A stable and accurate method for computing ID is proposed in~\cite{Cheng-Random}, similar
to the RRQR algorithm in~\cite{Gu-RRQR}.
We can combine RS with ID to get a more efficient low-rank approximation algorithm~\cite{Martinsson-PNAS07}.
For a given $n\times n$ matrix $B$, generate an $n\times (r+p)$ Gaussian random matrix $\Omega$ as above, and compute
the row sampling and column sampling matrices $Y=B\Omega$ and $Z = B^T \Omega$.
Then, use ID to determine the $r$ selected rows and columns of $B$ from $Y$ and $Z$,
\begin{equation}
[X^{row}, I^{row}] = \texttt{interpolative}(Y^T), \quad
[X^{col},J^{col}] = \texttt{interpolative}(Z^T),
\end{equation}
and $B$ can be approximated by
\[
B \approx X^{row} \cdotp B(I^{row},J^{col}) \cdotp (X^{col})^T.
\]

\subsection{Random HSS construction for Cauchy-like matrices}
\label{sec:rhss-cauchy}

The main idea is to apply the randomized ID to the row and column sampling matrices by
traversing the HSS tree level-by-level, from bottom to top.
To illustrate it, let $A$ be a matrix as defined in~\eqref{eq:Ev-Q}, $Y=A\Omega^{(1)}$ and $Z=A^T\Omega^{(2)}$
be the sampling matrices, where $\Omega^{(i)}$ is a Gaussian random matrix for $i=1,2$.
To construct an HSS matrix, we need to find the low-rank approximations of all HSS blocks, $H_i^{row}$ and
$H_i^{col}$. Recall that $H_i^{row}$ and $H_i^{col}$ are respectively the $i$-th HSS block row and column,
satisfying
\begin{equation}
A_{t_i\times \mathcal{I}} = H_i^{row}+D_i, \quad A_{\mathcal{I}\times t_i} = H_i^{col}+D_i.
\end{equation}

In this subsection, we show how to obtain the low-rank approximations from $Y$ and $Z$ by using the randomized ID method.
For a leaf node $i$, its compressed HSS block row and column are, respectively,
\[
\Phi_i=Y_i-D_i\Omega_i^{(1)}, \quad \Theta_i=Z_i-D_i^T\Omega_i^{(2)},
\]
where $(\star)_i$ means $(\star)(t_i,:)$ for $( \star )=Y,Z,\Omega^{(1)}$ and $\Omega^{(2)}$.
By applying the ID method to $\Phi_i$ and $\Theta_i$, we can easily obtain the low-rank approximations to $H_i^{row}$ and
$H_i^{col}$, respectively.

For a parent node, its compressed HSS blocks can be neatly obtained from those of its
children recursively, see section 4.1 in~\cite{rand-hss} and Algorithm~\ref{alg:randhss} below.
Then, its generators can be obtained similarly by applying ID to the compressed HSS blocks.

\begin{algorithm}
\label{alg:randhss}
(Random HSS construction for Cauchy-like matrices)
Given the generators of Cauchy-like matrix $A$, compute its HSS matrix approximation accurately.

First, use~\eqref{eq:app} to estimate the HSS rank $r$ of $A$ and generate two $N\times (r+p)$ Gaussian random matrices
$\Omega^{(1)}$ and $\Omega^{(2)}$. Then, compute $Y=A\Omega^{(1)}$ and $Z=A^T \Omega^{(2)}$.

\begin{description}

\item \texttt{do} $\ell=L,\cdots,1$
  \item[\quad] \texttt{for} node $i$ at level $\ell$
    \item[\qquad] \texttt{if} $i$ is a leaf node,
      \begin{enumerate}
        \item $D_i=A_{t_i,t_i}$;

        \item compute $\Phi_i=Y_i-D_i\Omega_i^{(1)}$, $\Theta_i=Z_i-D_i^T \Omega_i^{(2)}$;

        \item compute the ID of  $\Phi_i \approx U_i \Phi_i|_{\tilde{I}_i}$, $\Theta_i \approx V_i \Theta_i|_{\tilde{J}_i} $;

        \item compute $\widehat{Y}_i = V_i^T \Omega_i^{(1)}$, $\widehat{Z}_i = U_i^T \Omega_i^{(2)}$;
       \end{enumerate}

    \item[\qquad] \texttt{else}
      \begin{enumerate}

        \item store the generators $B_{i_1} = A(\tilde{I}_{i_1},\tilde{J}_{i_2})$, $B_{i_2} = A(\tilde{I}_{i_2},\tilde{J}_{i_1})$;

        \item compute $\Phi_i=\begin{bmatrix} \Phi_{i_1}|_{\tilde{I}_{i_1}}-B_{i_1}\widehat{Y}_{i_2} \\ \Phi_{i_2}|_{\tilde{I}_{i_2}}-B_{i_2}\widehat{Y}_{i_1} \end{bmatrix}$, \quad
          $\Theta_i=\begin{bmatrix} \Theta_{i_1}|_{\tilde{J}_{i_1}}-B_{i_2}^T\widehat{Z}_{i_2} \\ \Theta_{i_2}|_{\tilde{J}_{i_2}}-B_{i_1}^T\widehat{Z}_{i_1} \end{bmatrix}$;

        \item compute the ID of  $\Phi_i \approx U_i \Phi_i|_{\tilde{I}_i}$, \quad $\Theta_i \approx V_i \Theta_i|_{\tilde{J}_i}$;

        \item Compute $\widehat{Y}_i = V_i^T \begin{bmatrix} \widehat{Y}_{i_1} \\ \widehat{Y}_{i_2} \end{bmatrix}$,  \quad
                      $\widehat{Z}_i = U_i^T \begin{bmatrix} \widehat{Z}_{i_1} \\ \widehat{Z}_{i_2} \end{bmatrix}$;
       \end{enumerate}

    \item[\qquad] \texttt{end if}

   \item[\quad] \texttt{end for}

\item \texttt{end do}

\end{description}
For the root node $i$, store $B_{i_1} = A(\tilde{I}_{i_1},\tilde{J}_{i_2})$, $B_{i_2} = A(\tilde{I}_{i_2},\tilde{J}_{i_1})$.

\end{algorithm}


It can be verified that the complexity of Algorithm~\ref{alg:randhss} is $C_M+O(N r^2)$, where
$C_M$ is the cost of multiplying $A$ with (two) random matrices, and $r$ is the HSS rank of $A$.
In practice, we can let $\Omega^{(1)}=\Omega^{(2)}$.
FMM can be used to compute the sample matrices $Y$ and $Z$, which only costs $O((r+p)N \log N )$ flops.
For large matrices, FMM can be much faster than the plain matrix-matrix multiplications.
If using FMM, the complexity of RSHSS is $O(N r^2)$ flops,
see the reference~\cite{rand-hss}.
This HSS construction algorithm in theory can be faster than the algorithm proposed in~\cite{Shengguo-SIMAX2} which costs $O(N^2 r)$ flops.

Most time of RSHSS is taken to compute the sample matrices.
In the sequential case it takes about $80\%$ of the construction time and
about $30\%$ in the fully parallel case, refer to Table~\ref{tab:Ex1-speedups}.
In~\cite{Martinsson-Rev10}, it is proposed to use the subsampled random Fourier (SRFT) or Hadamard (SRHT) transforms to compute
the sample matrices.
We do not use this technique in RSHSS or ADC2, since the SRFT would introduce complex matrices, and
the SRHT requires the dimension of matrix $A$ to be powers of two.
Furthermore, the construction algorithm is usually much faster than the HSS matrix multiplication algorithm, see the
results in Table~\ref{tab:Ex1-speedups}.
Note that if the SRHT is applicable, the complexity of RSHSS is also about $O(Nr^2)$ flops.

Another issue we want to mention is the accuracy of RSHSS.
If the singular values of HSS blocks do not decay rapidly, the RS
method may lose a bit of accuracy. A power scheme was proposed to improve the quality of sample matrices in~\cite{Martinsson-Rev10},
for instance compute $Y=(AA^T)^qA \Omega$.
We find it is very difficult to incorporate this technique into
Algorithm~\ref{alg:randhss} and moreover, using the power scheme would require about $2q+1$ times as many
operations as Algorithm~\ref{alg:randhss}.
For accuracy, we choose a relatively large oversampling parameter $p$ and try to compute the ID of sampled matrices
as accurately as possible.
In practice, we let $\epsilon=1e$-$16$ in~\eqref{eq:app} to estimate the rank, and let the oversampling
parameter $p=10$.
The number of used random vectors is usually larger than the
HSS rank.
This strategy in practice is quite robust and it does not fail for any experiments during all our tests.
Note that RSHSS still has a risk of losing accuracy, for example, if $r+p$ is smaller than the HSS rank in some rare cases.

\section{Implementation details}
\label{sec:impl}

The ADC algorithm is consisted of three other algorithms:
\emph{the HSS construction} and \emph{HSS matrix multiplication algorithms},
and the standard DC algorithm.
Almost all modern CPUs have multiple cores, and we implement the ADC algorithms
in parallel to exploit the multicore architecture.
We use OpenMP to implement these algortihms.
This section introduces the parallel implementation details of these three algorithms.

\subsection{Parallel RSHSS algorithm}

As illustrated in Algorithm~\ref{alg:randhss} and section 4.1 of~\cite{rand-hss}, the computations for
different nodes at the same level can be performed simultaneously.
We can exploit the parallelism of the HSS tree, and
the computations for different nodes are done by different \emph{processes}.
Furthermore, the work for each node can also be done by \emph{multi-threads} by calling
a multithreaded BLAS library.

Recall that Algorithm~\ref{alg:randhss} computes three or four generators for each node,
$D_i$, $U_i$, $V_i$ and $B_i$. Note that parent nodes do not have
the generator $D_i$, and $U_i$ is computed from the row compression, and
$V_i$ is from the column compression. The generators $D_i$ and $B_i$ are submatrices
of the original matrix $A$, and are Cauchy-like.
Besides the number of flops,
the running time of algorithms is also determined by the amount of data movements.
To have good data locality,
we store the same type of generators for nodes at the same level continuously.
For example, we first store all the generators $U_i$ at level $\ell$,
then the generators $V_i$ and finally the generators $B_i$ at level $\ell$, for $\ell=L,\ldots,1$.
All the generators are stored continuously in one array, name it $A_H$, and the
generators $D_i$ are stored in the front part of $A_H$.
This form of storage is good for HSS matrix multiplications, see Algorithm~\ref{alg:phssmm-omp} below,
where the computations follow the HSS tree level by level, and
the generators of the same type are used one after the other.
For example, the computations of ~\eqref{eq:compt_X} use all the generators $D_i$ at the bottom level, and
so do the generators $V_i$.

Another point we want to mention is that the Cauchy-like matrices $D_i$ and $B_i$ are computed
from its generators respectively,
which are four vectors, see equation~\eqref{eq:cauchylike}.
We find that recomputing the entries of $D_i$ and $B_i$ is usually faster than subtracting them from
the original matrix $A$.

Our parallel version of RSHSS is similar to Algorithm~\ref{alg:randhss}.
The only difference is that the \texttt{do-loop} in Algorithm~\ref{alg:randhss} is replaced
by the following process after some computation details are ignored.

\begin{description}

\item \texttt{par\_for} leaf node $i$,
\item \hspace{0.5cm}  compute the Cauchy-like matrix $D_i$ via its generators and store it in $A_H$;
\item \texttt{end par\_for}

\item \texttt{do} $\ell=L,\cdots,1$
  \item \texttt{par\_for} node $i$ at level $\ell$, compute its generator $U_i$ from  $\Phi_i$ and
    store $U_i$ in $A_H$; \texttt{end par\_for}

  \item \texttt{par\_for} node $i$ at level $\ell$, compute its generator $V_i$ from  $\Theta_i$ and
    store $V_i$ in $A_H$; \texttt{end par\_for}

  \item \texttt{par\_for} node $i$ at level $\ell$, compute the Cauchy-like matrix $B_i$ and
    store it in $A_H$; \texttt{end par\_for}

\item \texttt{end do}

\end{description}

The abbreviation \texttt{par\_for} stands for `parallel for', which means the following computations can be done in parallel.
In practice we use~\eqref{eq:app} to estimate the HSS rank $r$,
based on the partition of the original matrix $A$, see Figure~\ref{fig:Ahss}.
The matrix $A$ is partitioned by letting all leaf nodes have roughly $m$ rows and columns.
Since $\{d_i\}$ and $\{\lambda_i\}$ are ordered increasingly, each partition of $Q$~\eqref{eq:Ev-Q} can also be seen as a
partition of interval $[\hat{a}, \hat{b}]$ which contains both $\{d_i\}$ and $\{\lambda_i\}$.
For the partition in Figure~\ref{fig:Ahss}, the interval $[\hat{a}, \hat{b}]$  is divided into four segments, and
the first $m_1$ entries of $\{d_i\}$ and $\{\lambda_i\}$ lie in the first segment of $[\hat{a}, \hat{b}]$,
the second $m_2$ entries lie in the second segment, and so on.
The rank estimated by~\eqref{eq:app} depends on the distance of two segments, which is defined in section~\ref{sec:lowrank}.
We use the distances of neighbouring segments to estimate rank, and choose the \emph{maximum} rank estimated
by~\eqref{eq:app} as the HSS rank.
For Figure~\ref{fig:Ahss}, there are \emph{three} pairs of neighbouring segments, and
the estimated ranks are of $H_1^{row}$, $H_3^{row}$ and $H_4^{row}$ respectively, and
the maximum of them is used as an estimate of HSS rank $r$.
If some eigenvalues are clustered, i.e., the distance between $I_i$ and $I_{i+1}$ is small,
the estimated rank by~\eqref{eq:app} may be too large to be useful.
We use the following tricks to get a more reasonable estimate of $r$.

\begin{enumerate}[(1)]
  \item If the distance between $I_i$ and $I_{i+1}$ is too small, we modify the partition of matrix $A$, i.e., move
    the boundary forward or backward to let the distance large.
    In our implementation, we modify the partition when the distance between $I_i$ and $I_{i+1}$ is less than $1e-10$,
    and the boundary is moved forward or backward by at most $k$(=5) rows and columns.
   \item If the computed rank by~\eqref{eq:app} is still too large, larger than 100, we fix the rank to be 100.
     (We find that HSS rank $r$ is rarely larger than 100 in the tridiagonal DC algorithm.)
\end{enumerate}

Note that these techniques are unfortunately lack of theoretical support, but
they make the rank estimation method more useful and robust.

\subsection{HSS matrix multiplication from the right}

After an HSS matrix is represented in its HSS form, there exist fast
algorithms for multiplying it with a vector in $O(Nr)$ flops
(see~\cite{ChandrasekaranDe06,Lyons-thesis}).
An HSS matrix multiplication algorithm has been introduced in~\cite{Lyons-thesis} for $H\times A$, where $H$ is an HSS matrix and $A$ is a general matrix.
For completeness, this subsection introduces the process of multiplying an HSS matrix with a general matrix
from right, i.e., compute $A\times H$. From Algorithm~\ref{alg:phssmm-omp} it is easy to see that the HSS matrix multiplication algorithms
are naturally parallelizable.

\begin{algorithm}[HSS matrix multiplication from right]
  \label{alg:phssmm-omp}
  Assume that the HSS tree $\mathcal{T}$ is a full binary tree and there are
  $L+1$ levels, the root is at level 0 and the leaf nodes are at level $L$.
  Let $j$ be the sibling of $i$.
  Let $X$ be a $P\times N$ matrix and partition the columns of $X$ as
  $X=[X_{i}]$, where $X_{i}=X(:,t_{i})$, $i\in LN$ is a leaf node.

  \begin{enumerate}[\quad (1)]

  \item \texttt{upsweep for} $G_i$
    \begin{itemize}
    \item
      \texttt{par\_for} $i$ at the bottom level, compute $G_i = X_i \cdotp U_i$; \texttt{end par\_for}

       \item \texttt{for} $\ell = L-1:-1:1$
         \begin{description}
         \item \texttt{par\_for} $i$ at level $\ell$, compute
           $G_i =\begin{bmatrix}  G_{i_1} & G_{i_2} \end{bmatrix}\cdotp  U_i$; \texttt{end par\_for}
         \end{description}

       \item \texttt{end for}
       \end{itemize}

     \item \texttt{downsweep for} $F_i$
       \begin{itemize}
       \item \texttt{par\_for} $i$ at the second top level, compute
        $F_i=G_j\cdotp B_j, \begin{bmatrix} F_{i_1} \\ F_{i_2}\end{bmatrix} = F_i \cdotp V_i^T;$
         \texttt{end par\_for}

       \item \texttt{for} $\ell = 2:L-1$
         \begin{description}
         \item \texttt{par\_for} $i$ at level $\ell$, compute $F_i=G_jB_j+F_i,$
           $\begin{bmatrix} F_{i_1} \\ F_{i_2} \end{bmatrix}=F_i\cdotp V_i^T;$  \texttt{end par\_for}
         \end{description}

       \item \texttt{end for}
       \end{itemize}

     \item compute $X$
       \begin{itemize}
       \item \texttt{par\_for} $i$ at the bottom level, compute
         \begin{equation}
           \label{eq:compt_X}
           X_i = X_iD_i+F_iV_i^T;
         \end{equation}

       \item \texttt{end par\_for}
       \end{itemize}
  \end{enumerate}
\end{algorithm}

All the computations for the nodes at the same level are independent of each other.
Furthermore, almost all the operations are matrix-matrix multiplications and we can take advantage of
the highly optimized routine \texttt{DGEMM} in MKL.
We explore both the parallelism in the HSS tree and the parallelism
from the blas operations by using MKL.
Table~\ref{tab:hss_parallel} shows the speedups of Algorithm~\ref{alg:phssmm-omp} when
only exploiting the parallelism in the HSS tree.
The dimension of the HSS matrix is 10000, which is defined in the same way
as the matrix $Q$ in Example 1, and we multiple it with a $10000\times 10000$
random matrix via Algorithm~\ref{alg:phssmm-omp}.
The times cost by Algorithm~\ref{alg:phssmm-omp} are presented in the third row of Table~\ref{tab:hss_parallel},
and the compiled codes are linked to a sequential BLAS library.
The results in Table~\ref{tab:hss_parallel} show that the scalability of Algorithm~\ref{alg:phssmm-omp} is good.
Some more numerical results are included in Example 2 in section~\ref{sec:numer-tedc}.

\begin{table}[ptbh]
\caption{The parallelism of Algorithm~\ref{alg:phssmm-omp} introduced by the HSS tree structure}%
\label{tab:hss_parallel}
\begin{center}%
\begin{tabular}
[c]{|c|c|c|c|c|c|c|c|c|c|} \hline
\multirow{2}{*}{} & \multicolumn{9}{c|}{Threads}  \\ \cline{2-10}
                  & $1$ & $3$  & $5$  & $7$  & $9$  & $11$ & $13$ & $15$ & $16$  \\ \hline \hline
time(s)           &  11.34 & 4.06 & 2.63 & 2.07 & 1.73 & 1.51 & 1.40 & 1.31 & 1.20   \\ \hline
speedups          &  1.00  & 2.79 & 4.31 & 5.48 & 6.55 & 7.51 & 8.10 & 8.66 & 9.45  \\ \hline
\end{tabular}
\end{center}
\end{table}

\subsection{Accelerate the process of DC algorithm}

The LAPACK routine \texttt{dstevd} implements a divide-and-conquer algorithm for
symmetric tridiagonal matrices.
It computes the eigenvalues and eigenvectors
explicitly by calling \texttt{dlaed0}.
The routine \texttt{dlaed0} solves each subproblem in a divide-and-conquer way.
\texttt{dlaed1} called by \texttt{dlaed0} computes the eigendecomposition of the merged subproblem,
and it calls \texttt{dlaed2} to deflate a diagonal matrix with rank-one modification and calls \texttt{dlaed3}
to update the eigenvector matrix via matrix-matrix multiplications.

Our implementation has the same structure as LAPACK.
We add the HSS techniques in the routine \texttt{dlaed1} and rename it \texttt{mdlaed1}.
When the size of the deflated matrix is small, it calls \texttt{dlaed3} as usual.
Otherwise, it calls \texttt{mdlaed3} to compute the eigenvalues and
update the eigenvectors.
In our implementation, we use the HSS matrix techniques when the size of the deflated matrix
is larger than 2000. The routine \texttt{mdlaed3} is similar to
\texttt{dlaed3}, and it computes the
eigenvalues $\{\lambda_i\}$, the recomputed vector $\{\hat{u}_i\}$,
$\gamma_i=\lambda_i-d_i$ (the distance between $\lambda_i$ and $d_i$),
and $\mu_i=d_{i+1}-\lambda_i$ (the distance between $\lambda_i$ and $d_{i+1}$), for $i=1,\ldots,n$.
The secular equation in \texttt{mdlaed3} is solved in parallel by calling \texttt{dlaed4}.
Then the eigenvector matrix of the diagonal matrix with rank-one modification is approximated
by an HSS matrix, and the eigenvectors of the original matrix $T$ are updated via
fast HSS matrix multiplications~\cite{Lyons-thesis}, see also Algorithm~\ref{alg:phssmm-omp}.
Using the HSS matrix techniques to update the eigenvectors saves a lot of flops since
the complexity is reduced from $O(N^3)$ to $O(N^2r)$.

\begin{figure}[ptbh]
\centering
\subfigure[Comparion of flops estimated by Vtune]{
\includegraphics[width=2.5in,height=2.0in]{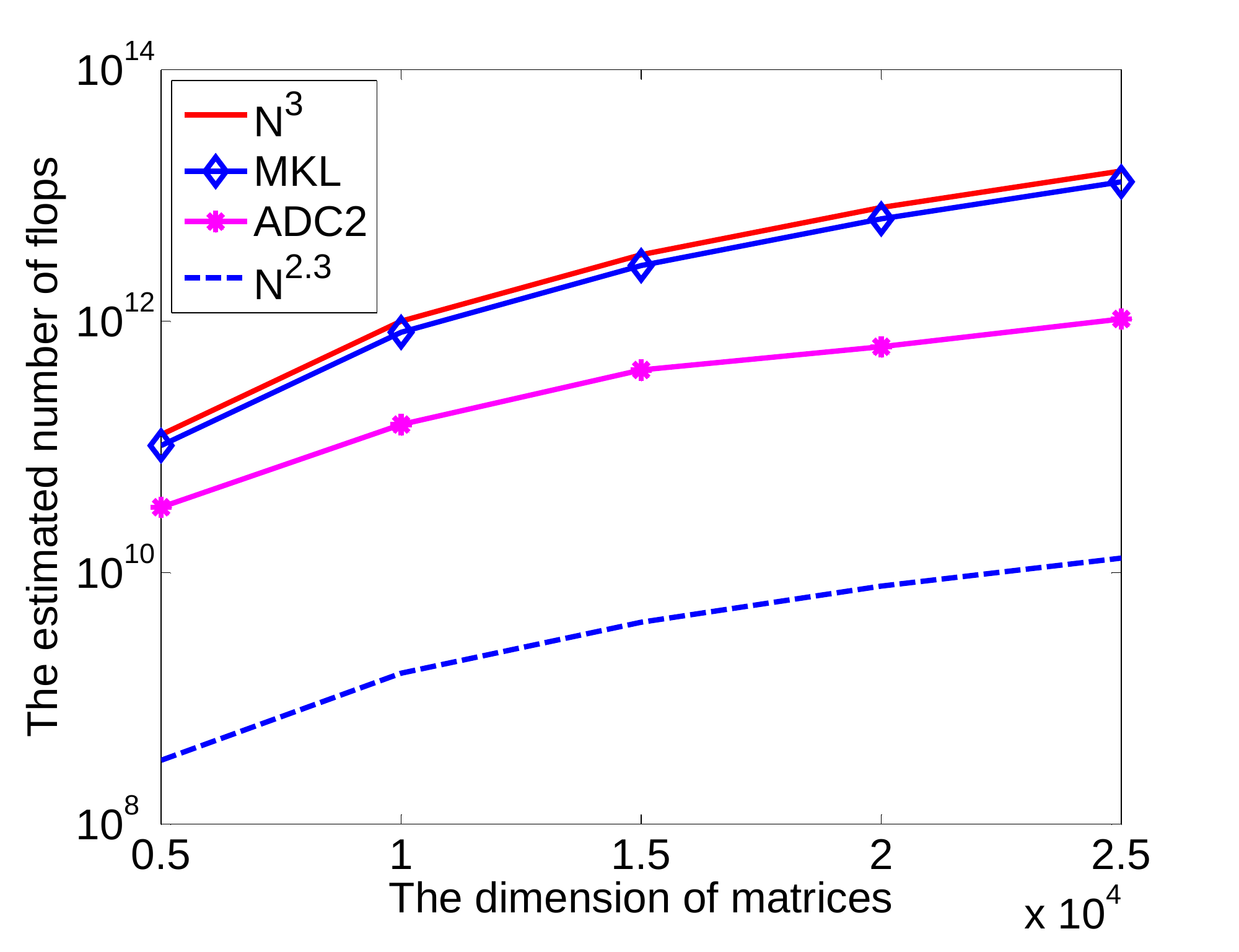}
\label{fig:flops}}
\subfigure[The divide and conquer tree]{
\includegraphics[width=2.0in,height=1.6in]{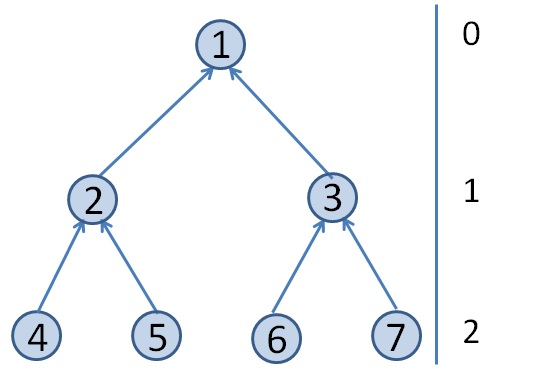}
\label{fig:comput_tree}}
\caption{The comparison of flops and the DC tree}%
\label{fig:flops_res}
\end{figure}

The divide-and-conquer algorithm is also organized in a binary tree structure.
Figure~\ref{fig:comput_tree} shows the tree structure of DC algorithm.
A big problem is recursivley splitted into two small problems, and two small problems
are merged together into a big one.
The subproblems at the same level of DC tree can be solved
in parallel.
The subproblems at the bottom level are solved by using the QR algorithm in parallel, calling
the LAPACK routine \texttt{dlasdq} in our implementation.
For the problems at the other levels, we
had tried to use the nested parallelism of OpenMP to exploit both the parallelism of DC tree and
BLAS operations, but it did not give us any speedup increases.
Furthermore, if using nested parallel computing, 
each thread would require a lot of private memory to store the intermediate eigenvectors, and the
memory cost would be greatly increased.
Therefore, the problems above the bottom level of DC tree are solved sequentially and
we only exploit the parallelism of BLAS operations and HSS techniques.

As is well known, solving the secular equations costs about $O(N^2)$ operations.
We further parallelize the process of solving the secular equations, which is
inspired by the work~\cite{DC_quark}.
We simply add OMP PARALLEL DO directives in \texttt{mdlaed3} when calling \texttt{dlaed4}.
Our parallel implementation is simpler than that in~\cite{DC_quark}.
In this paper we use OpenMP and follow the fork-join model.
While, it followed a task-flow model in~\cite{DC_quark} and used a dynamic runtime system to schedule the tasks.
Comparing the numerical results in~\cite{DC_quark} with those in the next section, we can see that the
rank-structured matrix techniques is good for the case that there are few deflations and that
the task-flow model used in~\cite{DC_quark} is good for the case that there are a lot of
deflations.
The advantage of the task-flow model is that it introduces a huge level of parallelism
through a fine task granularity, and that the tasks are scheduled by a runtime system, 
some synchronization barriers are removed.
As the algorithm in~\cite{DC_quark} still uses plain matrix-matrix multiplications to update the eigenvectors,
when the deflations are few, its the advantage decreases as
the dimensions of matrices increase.
Therefore, a good research direction is to combine the rank-structured matrix techniques with the task-flow model.

\section{Numerical results}
\label{sec:numer-tedc}

All the results are obtained on a server with 128GB memory and an Intel(R) Xeon(R) CPU E5-2670, which has
two sockets, 8 cores per socket, and 16 cores in total.
The codes are written in Fortran 90.
For compilation we used Intel fortran compiler (\texttt{ifort}) and the optimization flag
\texttt{-O2 -openmp}, and then linked the codes to Intel MKL (composer\_xe\_2013.0.079).

{\bf Example 2.} In this example, we use the matrix defined in Example 1 to show the
scalability of the HSS construction and matrix multiplication algorithms when implemented in parallel by using multi-threading.
The dimension of this matrix is 10000.
The row dimensions of HSS blocks for the leaf nodes are around $200$. 
The scalability of the HSS construction algorithm based on SRRSC and RSHSS are tested,
and the results are shown in Table~\ref{tab:Ex1-speedups}.
The results for Algorithm~\ref{alg:phssmm-omp} are also included.
The elapsed times of HSS constructions are shown in the rows denoted by \texttt{Const}, and
the times of HSS multiplications are included in those denoted by \texttt{Mult}.
The row denoted by \texttt{DGEMM} in Table~\ref{tab:Ex1-speedups} shows the times of computing the
sample matrices $Y$ and $Z$.

The results in Table~\ref{tab:Ex1-speedups} are obtained by letting \texttt{OMP\_NUM\_THREADS} and \texttt{MKL\_NUM\_THREADS}
equal to $1, 3, 5, \cdots, 15$ and $16$.
From the results we can see that the HSS construction algorithm is usually faster than the HSS multiplication algorithm.
For the RSHSS algorithm, we let $p$ equal to $10$ and the estimated rank by~\eqref{eq:app} is 79 which is
larger than 57, the HSS rank computed by SRRSC.
Most ranks of the HSS blocks are around 40.
The HSS matrix multiplications for RSHSS is slower than those for SRRSC,
since the ranks of HSS blocks computed by RSHSS are usually larger than those computed by SRRSC.
From the results in Table~\ref{tab:Ex1-speedups}, we can see that our parallel implementation
achieves good speedups.

\begin{table}[ptbh]
\caption{The execution time of HSS algorithms in seconds}%
\label{tab:Ex1-speedups}
\begin{center}%
\begin{tabular}
[c]{|c|c|c|c|c|c|c|c|c|c|c|}\hline
\multicolumn{2}{|c|}{ \multirow{2}{*}{Method} } & \multicolumn{9}{c|}{Threads}  \\ \cline{3-11}
\multicolumn{2}{|c|}{}         &   $1$ & $3$  & $5$  & $7$  & $9$  & $11$ & $13$ & $15$ & $16$  \\ \hline \hline
\multirow{2}{*}{SRRSC} & Const &  2.29 & 0.84 & 0.52 & 0.41 & 0.34 & 0.33 & 0.32 & 0.31 & 0.31   \\ \cline{2-11}
                       & Mult  &  4.75 & 1.87 & 1.32 & 1.08 & 0.96 & 0.97 & 0.89  & 0.85 & 0.85   \\ \hline \hline
\multirow{3}{*}{RSHSS} & DGEMM &  2.24 & 0.76 & 0.46 & 0.34 & 0.30 & 0.25 & 0.21 & 0.19 & 0.16  \\ \cline{2-11}
                       & Const &  2.85 & 1.07 & 0.75 & 0.63 & 0.57 & 0.53 & 0.48 & 0.45 & 0.42  \\ \cline{2-11}
                       & Mult &  8.04 & 2.98  & 1.94 & 1.54 & 1.29 & 1.14 & 1.09 & 1.02 & 0.95  \\ \hline
\end{tabular}
\end{center}
\end{table}

%

{\bf Example 3.}
For several classes of matrices~\cite{A880}, few or no eigenvalues are deflated
in the DC algorithm. Some of such matrices include
the Clement-type, Legendre-type, Laguerre-type, Hermite-type and Toeplitz-type matrices, which are defined as follows.
We use these matrices to show the performance of the ADC algorithms.

The Clement-type matrix~\cite{A880} is given by
\begin{equation}
  \label{eq:Clement-Tri}
  T=\text{tridiag} \small
  \begin{pmatrix}
    &\sqrt{n} & & \sqrt{2(n-1)} & & \sqrt{(n-1)2} & & \sqrt{n} & \\
    0 & & 0 & & \ldots & & 0 & & 0 \\
    &\sqrt{n} & & \sqrt{2(n-1)} & & \sqrt{(n-1)2} & & \sqrt{n} & \\
  \end{pmatrix},
\end{equation}
where the off-diagonal entries are $\sqrt{i(n+1-i)}, i=1,\ldots,n$.

The Legendre-type matrix is defined as~\cite{A880,Handbook-Func},
\begin{equation}
  \label{eq:Legendre-Tri}
  T=\text{tridiag} \small
  \begin{pmatrix}
    &2/\sqrt{3\cdotp 5} & & 3/\sqrt{5\cdotp 7} & & n/\sqrt{(2n-1)(2n+1)} & \\
    0 & & 0 & & \ldots & & 0 \\
    &2/\sqrt{3\cdotp 5} & & 3/\sqrt{5\cdotp 7} & & n/\sqrt{(2n-1)(2n+1)} & \\
  \end{pmatrix},
\end{equation}
where the off-diagonal entries are $i/\sqrt{(2i-1)(2i+1)}, i=2,\ldots,n$.

The Laguerre-type matrix is defined as~\cite{A880},
\begin{equation}
  \label{eq:Laguerre-Tri}
  T=\text{tridiag} \small
  \begin{pmatrix}
    &2 & & 3 & & n-1 & & n & \\
    3 & & 5 & & \ldots & & 2n-1 & & 2n+1 \\
    &2 & & 3 & & n-1 & & n & \\
  \end{pmatrix}.
\end{equation}

The Hermite-type matrix is given as~\cite{A880},
\begin{equation}
  \label{eq:Hermite-Tri}
  T=\text{tridiag} \small
  \begin{pmatrix}
    &\sqrt{1} & & \sqrt{2} & & \sqrt{n-2} & & \sqrt{n-1} & \\
    0 & & 0 & & \ldots & & 0 & & 0 \\
    &\sqrt{1} & & \sqrt{2} & & \sqrt{n-2} & & \sqrt{n-1} & \\
  \end{pmatrix}.
\end{equation}

\begin{figure}[ptbh]
\centering
\subfigure[ADC1]{
\includegraphics[width=2.5in,height=2.0in]{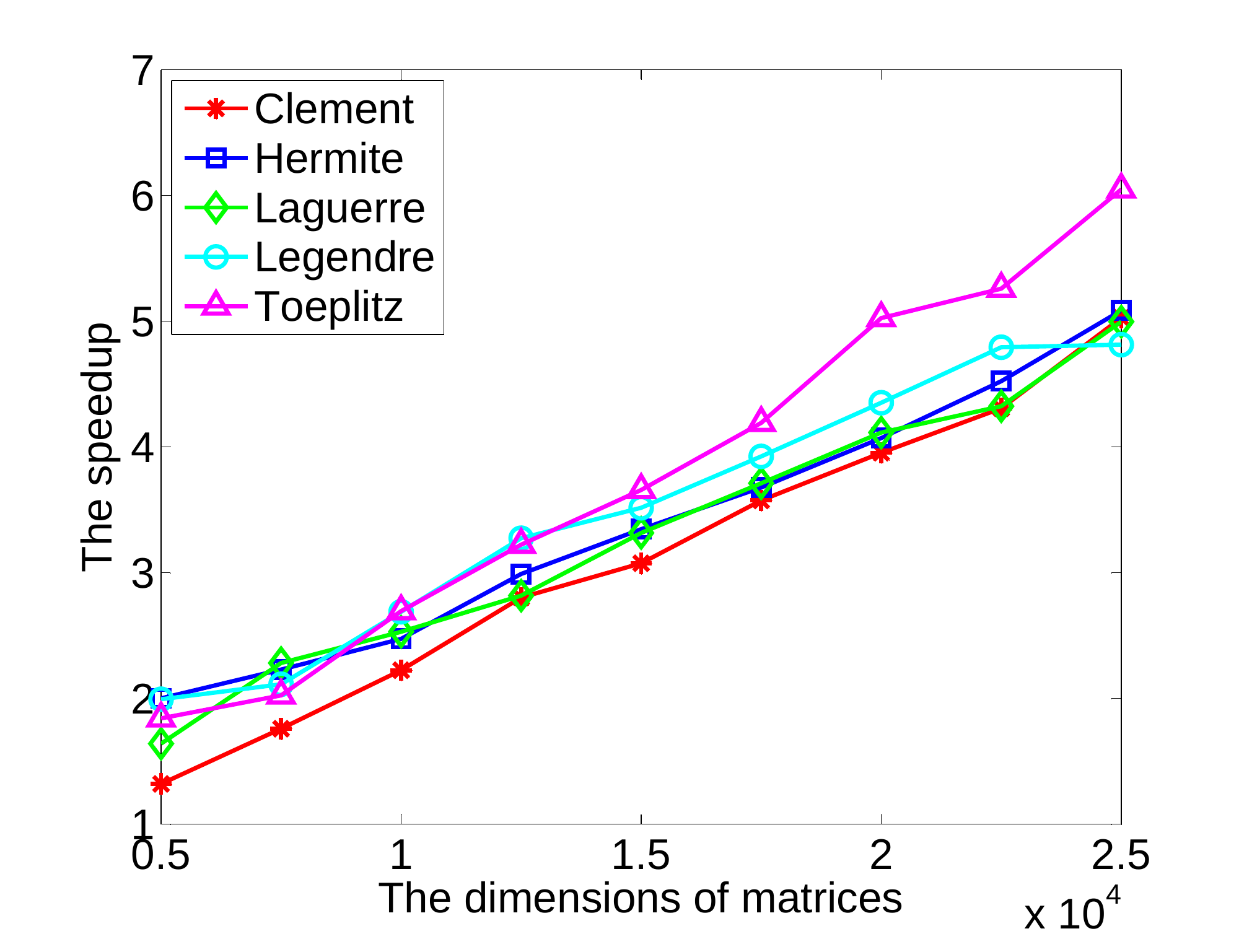}
\label{fig:spd-adc1}}
\subfigure[ADC2]{
\includegraphics[width=2.5in,height=2.0in]{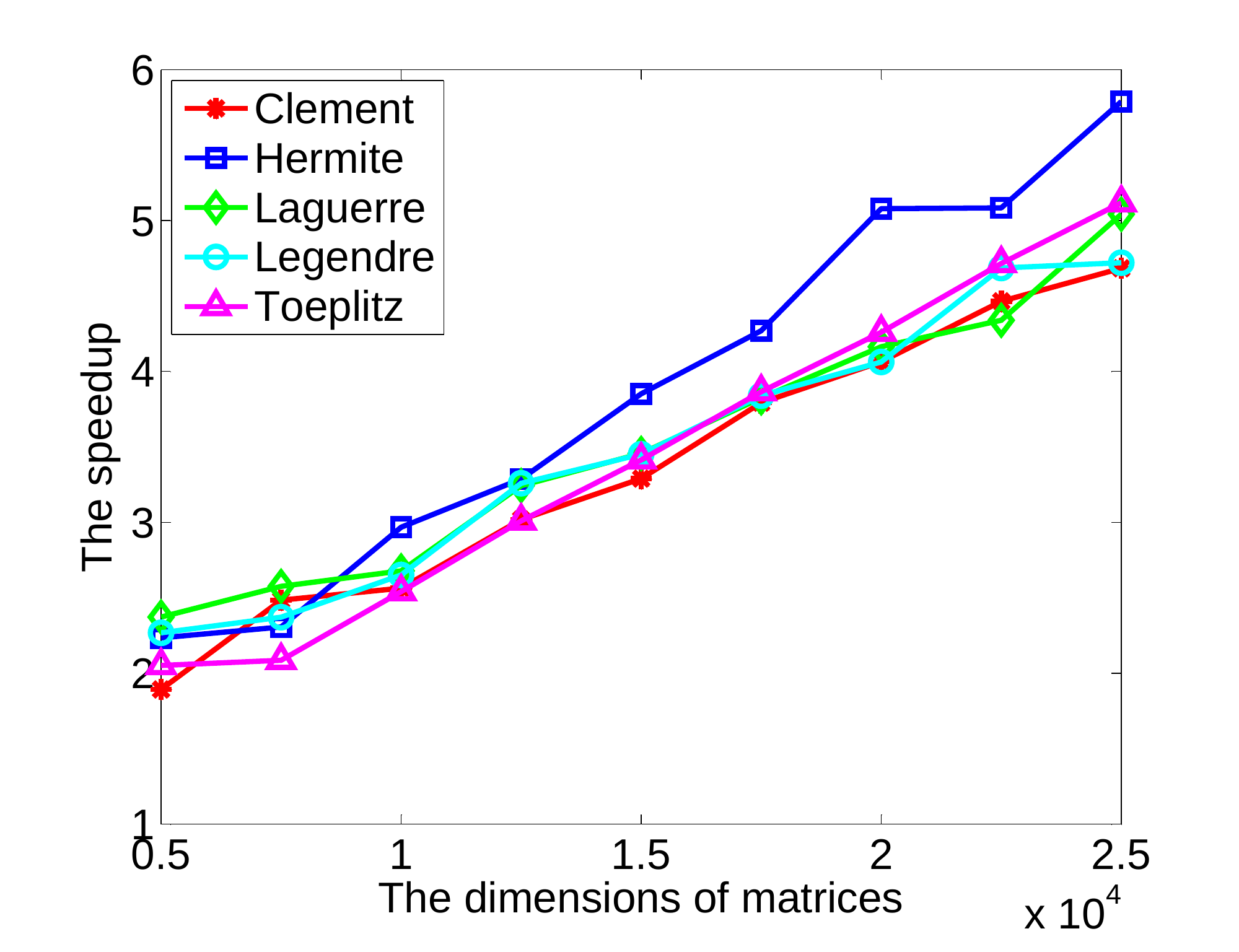}
\label{fig:spd-adc2}}
\caption{The weak scalability of ADC1 and ADC2}%
\label{fig:speedup}
\end{figure}

The Toeplitz-type matrix is a symmetric tridiagonal matrix with diagonals 2 and
off-diagonal entries 1.
In this example, we compare ADCs with DC in Intel MKL with
\texttt{OMP\_NUM\_THREADS}=16.
The speedups of ADC1 and ADC2 over DC are similar, and
the results are respectively reported in Table~\ref{tab:Ex2-srrsc} and Table~\ref{tab:Ex2-rhss},
see also Figure~\ref{fig:speedup}.
Since ADCs require fewer flops than the standard DC, they can achieve even better speedups for
larger matrices.
For example, when the dimension of Toeplitz-type matrix increases from $25k$ to
$40k$, the speedup of ADC2 over DC increases from $5.79$ to $8.06$.
During our experiments, HSS techniques are only used when the size of the current secular equation is
larger than 2000, which is a parameter depending on the computer architecture,
compiler and the optimized BLAS library.
The row dimensions of the HSS blocks for the leaf nodes are also around 200, and
16 threads are used.

\begin{table}[ptbh]
\caption{The speedups of {ADC1} compared with Intel MKL ($k$ denotes one thousand)}
\label{tab:Ex2-srrsc}
\begin{center}%
\begin{tabular}
[c]{|c|ccccccccc|}\hline
\multirow{2}{*}{Matrix}  & \multicolumn{9}{c|}{Dim} \\ \cline{2-10}
  & $5k$ & $7.5k$& $10k$ & $12.5k$ & $15k$ & $17.5k$ & $20k$ &$22.5k$ & $25k$  \\ \hline \hline
Clement  & 1.32x & 1.76x & 2.22x & 2.80x & 3.10x & 3.57x & 3.95x & 4.31x & 5.03x   \\
Legendre & 1.99x & 2.11x & 2.69x & 3.27x & 3.52x & 3.92x & 4.35x & 4.79x & 4.81x  \\
Laguerre & 1.64x & 2.28x & 2.53x & 2.82x & 3.32x & 3.71x & 4.11x & 4.32x & 5.00x \\
Hermite  & 2.00x & 2.23x & 2.47x & 2.99x & 3.35x & 3.68x & 4.07x & 4.52x & 5.08x  \\
Toeplitz & 1.84x & 2.02x & 2.69x & 3.22x & 3.66x & 4.19x & 5.02x & 5.26x & 6.05x \\ \hline
\end{tabular}
\end{center}
\end{table}

\begin{figure}[ptbh]
\centering
\subfigure[Max. error]{
\includegraphics[width=2.5in,height=2.0in]{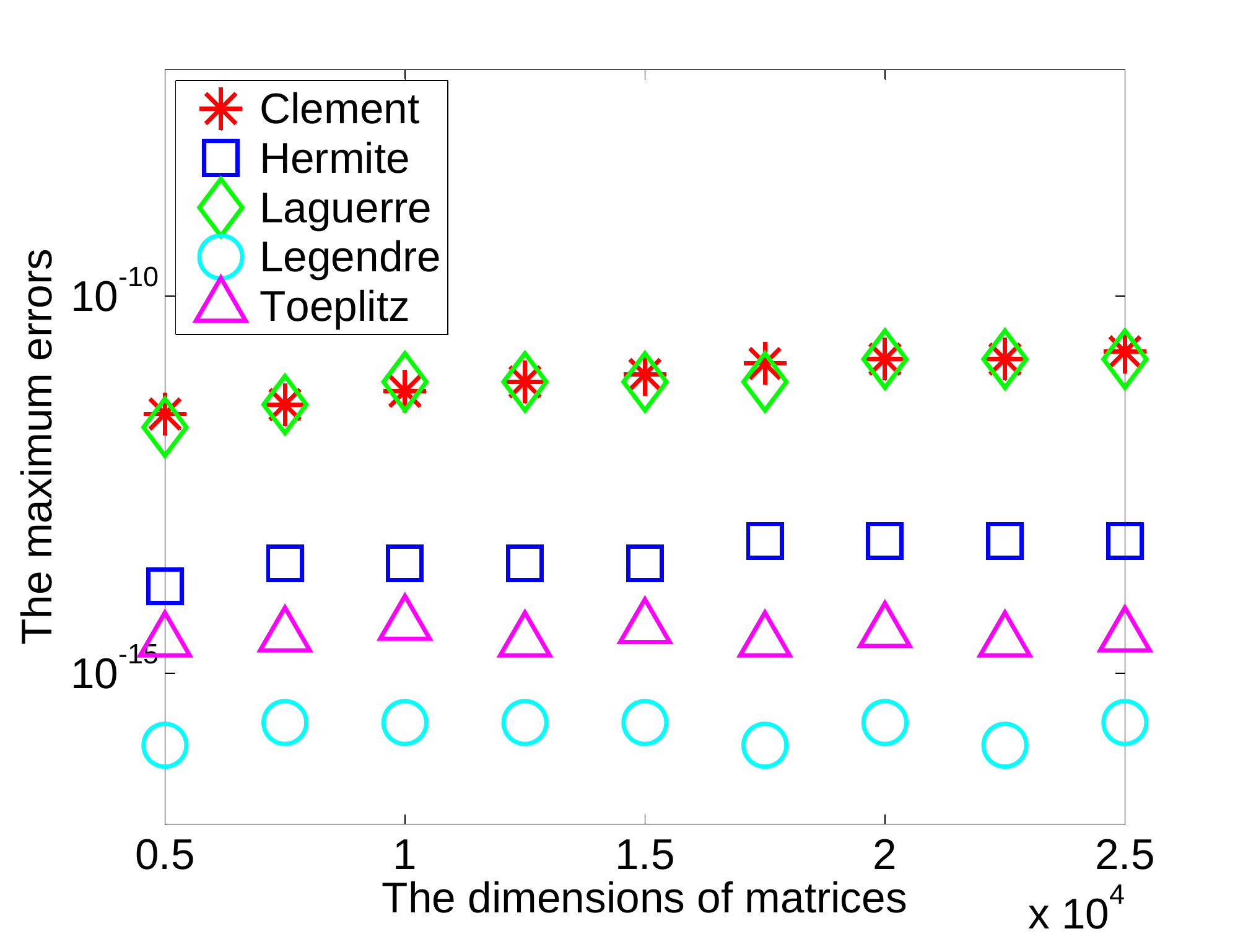}
\label{fig:maxerr}}
\subfigure[Max. Relative error]{
\includegraphics[width=2.5in,height=2.0in]{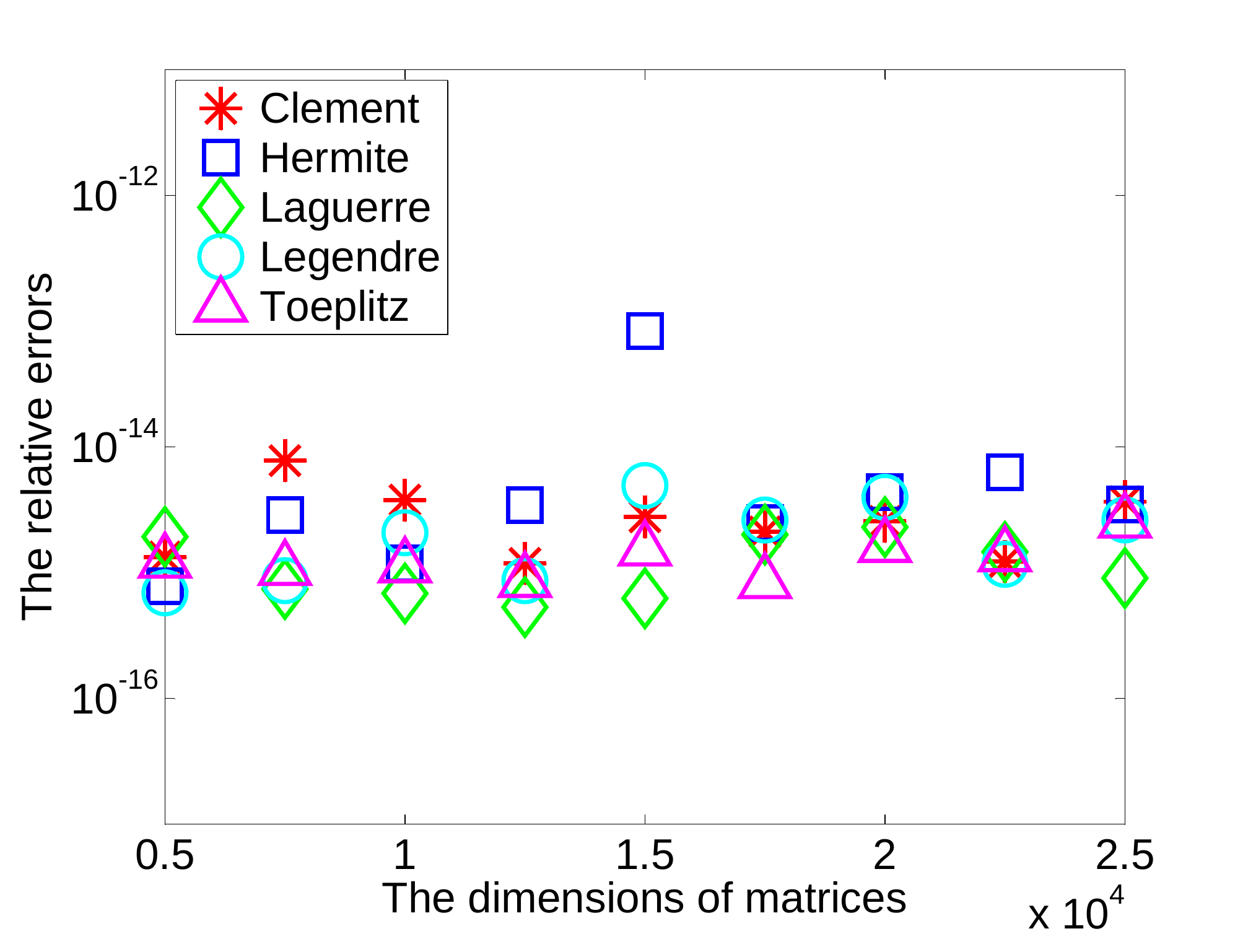}
\label{fig:relerr}}
\caption{Errors of the eigenvalues computed by ADC1 compared with those by MKL}%
\label{fig:Ex1-errors}%
\end{figure}

The DC algorithm is relatively complex, there are deflations and the secular
equations are solved via iterative methods,
and it is difficult to count the total number of flops cost by DC or ADCs by hand.
We use some tools based on event based sampling (EBS) technology to estimate the floating point operations.
PAPI~\cite{papi_web} and Intel Vtune Amplifier XE~\cite{vtune_web} are popular performance analysis tools, which can make reasonable estimates.
Figure~\ref{fig:flops} shows the comparisions of flops costed by ADC2 and DC in MKL.
These results are for the Toeplitz-type matrices with different dimensions.
From it we can see that DC requires nearly $O(N^3)$ flops, while
ADC2 requires much fewer flops.
Figure~\ref{fig:maxerr} and~\ref{fig:relerr} shows the maximum errors and maximum relative errors of the eigenvalues
computed by ADC1 compared with those by DC, respectively.
From the results, we can see that the computed eigenvalues by ADC1 nearly have the same accuracy as those computed by MKL.
Note that the results for relative error are included here but the DC algorithms in general are
not guaranteed to have high relative accuracy.

\begin{table}[ptbh]
\caption{The speedups of {ADC2} compared with Intel MKL ($k$ denotes one thousand)}
\label{tab:Ex2-rhss}
\begin{center}%
\begin{tabular}
[c]{|c|ccccccccc|}\hline
\multirow{2}{*}{Matrix}  & \multicolumn{9}{c|}{Dim} \\ \cline{2-10}
  & $5k$ & $7.5k$& $10k$ & $12.5k$ & $15k$ & $17.5k$ & $20k$ &$22.5k$ & $25k$  \\ \hline \hline
Clement  & 2.05x & 2.08x & 2.54x & 3.01x & 3.41x & 3.87x & 4.26x & 4.71x & 5.12x   \\
Legendre & 2.27x & 2.37x & 2.65x & 3.25x & 3.45x & 3.84x & 4.06x & 4.68x & 4.72x  \\
Laguerre & 1.89x & 2.48x & 2.56x & 3.01x & 3.28x & 3.79x & 4.06x & 4.46x & 4.68x \\
Hermite  & 2.37x & 2.58x & 2.68x & 3.33x & 3.46x & 3.82x & 4.16x & 4.33x & 5.04x  \\
Toeplitz & 2.23x & 2.31x & 2.97x & 3.29x & 3.85x & 4.27x & 5.07x & 5.08x & 5.79x \\ \hline
\end{tabular}
\end{center}
\end{table}

\begin{figure}[ptbh]
\centering
\subfigure[Orthogonality $\frac{\|I-QQ^T\|}{N}$]{
\includegraphics[width=2.5in,height=2.0in]{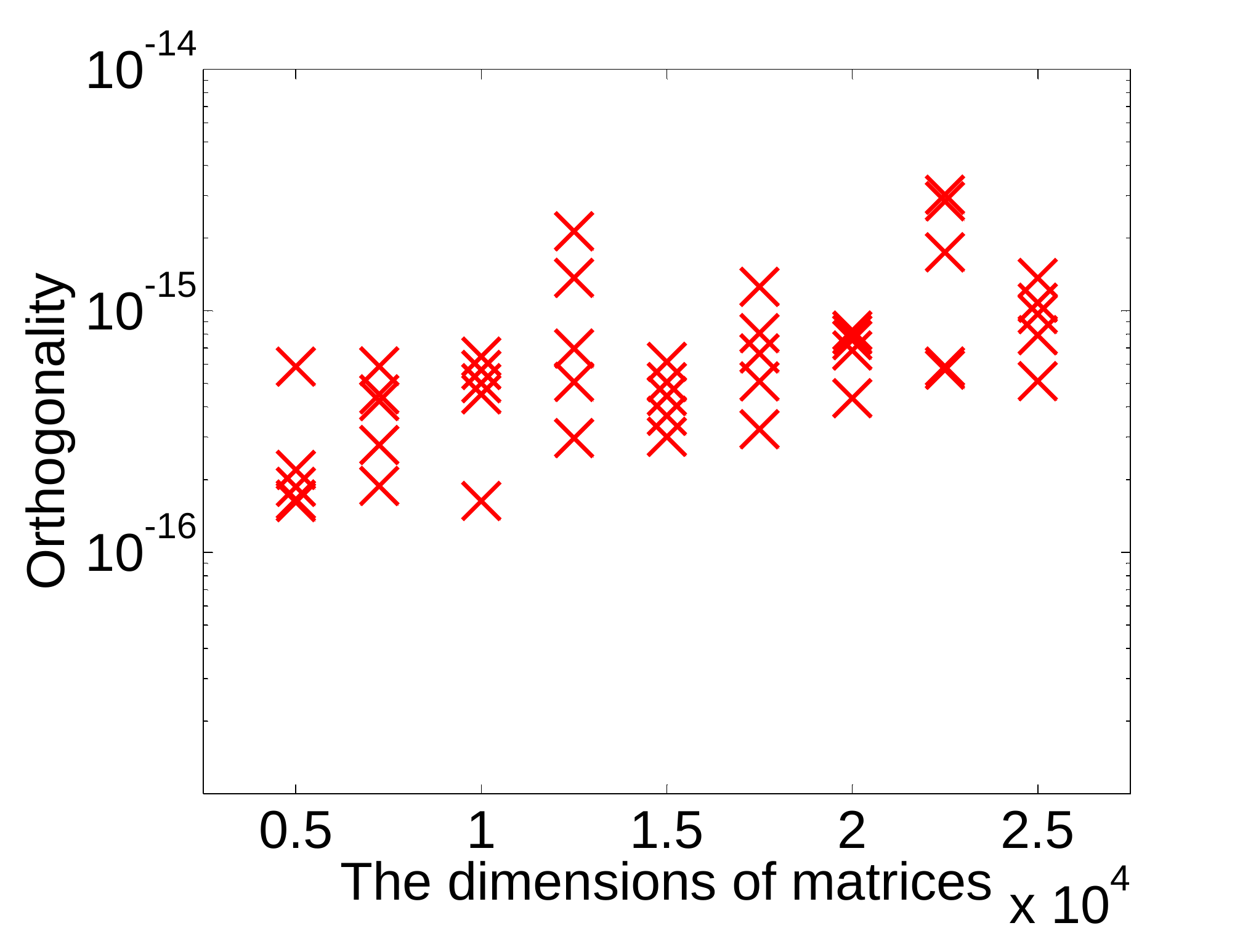}
\label{fig:orth-adc2}}
\subfigure[The backward error $\frac{\|T-Q\Lambda Q^T\|}{\|T\|\times N}$]{
\includegraphics[width=2.5in,height=2.0in]{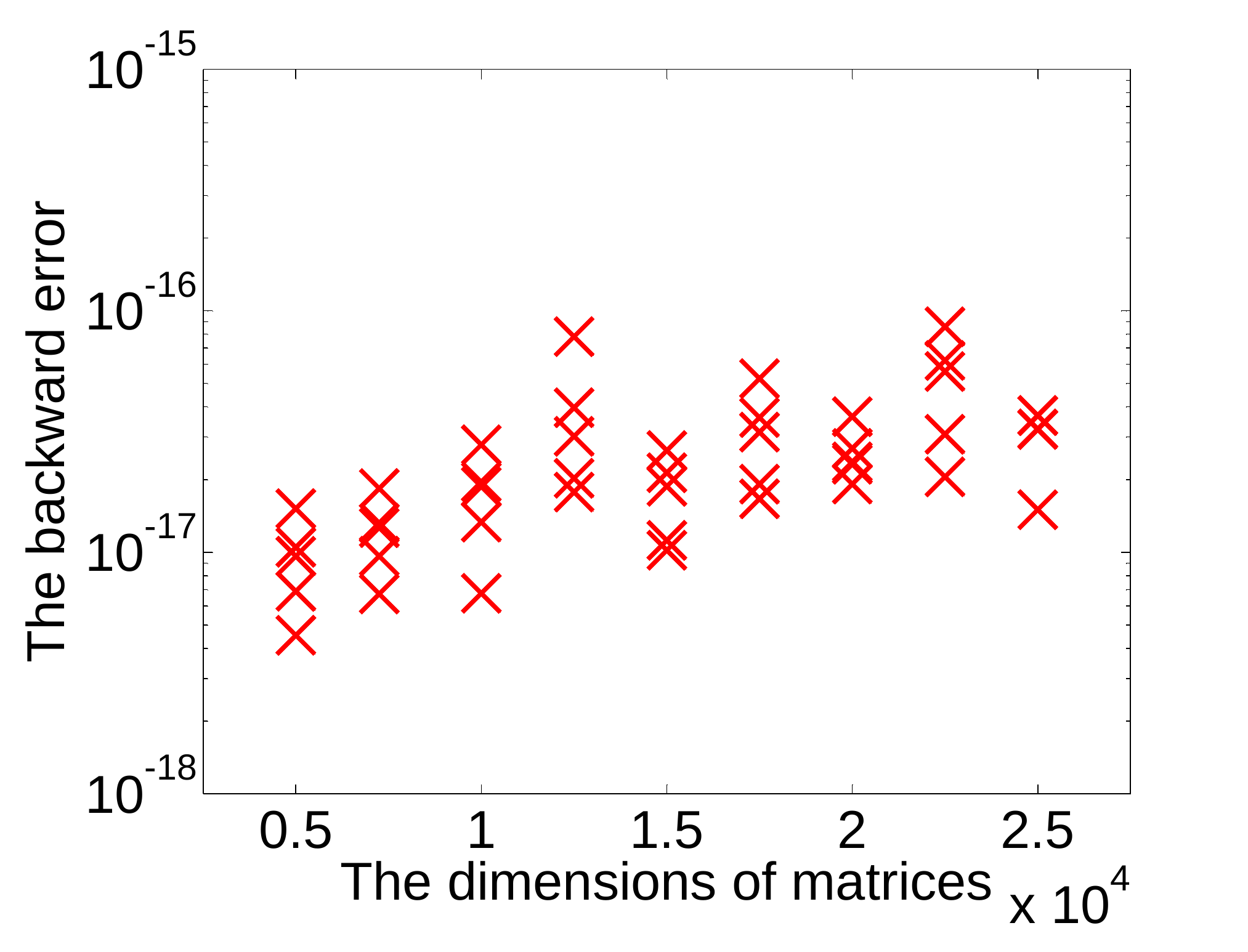}
\label{fig:back-adc2}}
\caption{The correctness of the eigenvectors computed by ADC2}%
\label{fig:orth-check}
\end{figure}

The results for the orthogonality of the computed eigenvectors are shown in Figure~\ref{fig:orth-adc2}, which are defined as $\frac{\|I-QQ^T\|}{N}$.
Figure~\ref{fig:back-adc2} shows the results for the backward error of ADC2, computed as $\frac{\|T-Q\Sigma Q^T\|}{\|T\|\times N}$.
While, ADC2 is a little less accurate than ADC1 but ADC2 can also be used reliably for applications,
the orthogonality of the computed eigenvectors by ADC2 is about 1$e$-12 and the maximum error
of the computed eigenvalues by ADC2 compared with those by Intel MKL is about 1$e$-14.
One advantage of ADC2 over ADC1 is that it requires fewer flops when FMM or SRHT is applicable,
which will be done in the future work.

\begin{figure}[ptbh]
\centering
\subfigure[The speedups of ADC2 compared with MKL]{
\includegraphics[width=2.5in,height=2.0in]{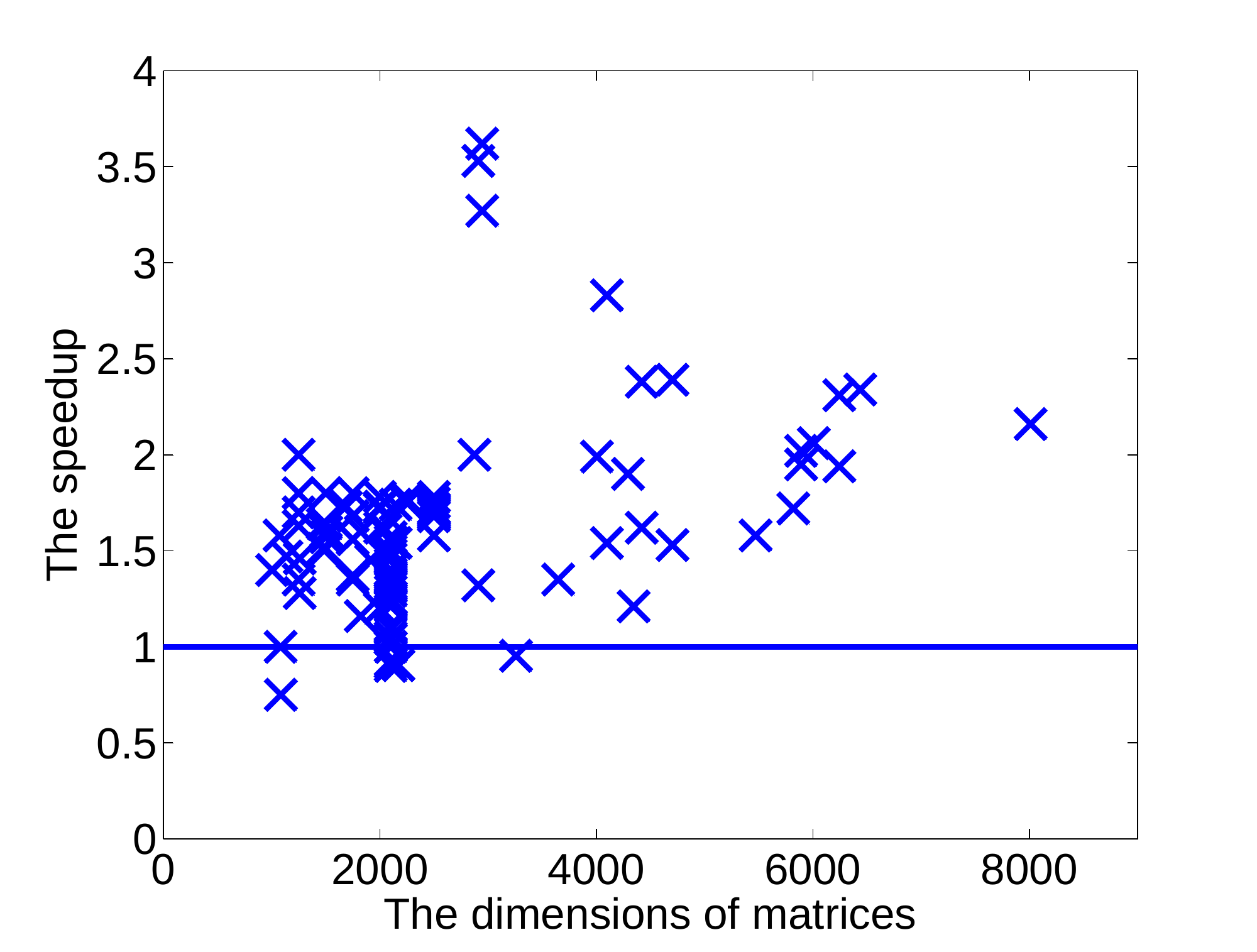}
\label{fig:orth-lapack}}
\subfigure[The relative errors of computed eigenvalues]{
\includegraphics[width=2.5in,height=2.0in]{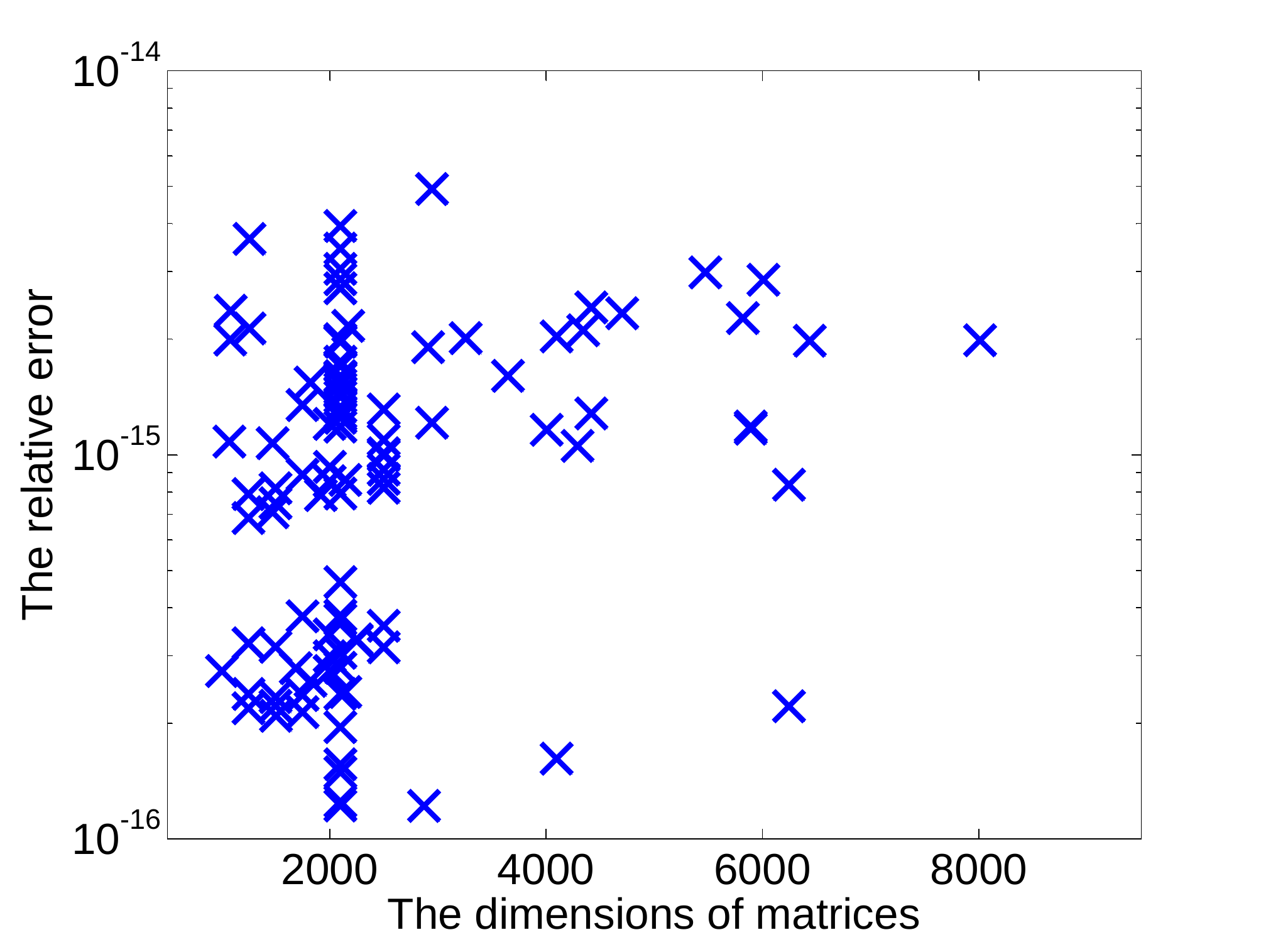}
\label{fig:spd-lapack}}
\caption{The results for matrices in \texttt{stetester}}
\label{fig:lapack-tester}
\end{figure}

We further use all the matrices in the LAPACK \texttt{stetester}~\cite{A880} with dimensions larger than
1000 to test ADC2.
Figure~\ref{fig:lapack-tester} shows the speedups of ADC2 over DC in MKL and the relative errors
of the eigenvalues computed by ADC2 compared with those by MKL.
The results show that for almost all matrices ADC2 is faster than the DC implmentation in MKL
and that the computed eigenvalues are highly accurate compared with
those computed by DC in MKL.
The experiments are done by letting \texttt{OMP\_NUM\_THREADS} and \texttt{MKL\_NUM\_THREADS}
equal to $16$.
For some rare matrices ADC2 is a little slower than \texttt{dstevd} in MKL but never slower by
more than $0.8$e-$02$ seconds.
Note that the HSS techniques are only used when the size of secular equation
is larger than 2000.
For the matrices with dimensions from 1000 to 2000, the speedups of ADC2
over DC in MKL is due to that ADC2 computes the bottom subproblems of the
DC tree and the secular equations in parallel.

\section{Conclusions}

In this paper, two accelerated tridiagonal DC algorithms are proposed by
using the HSS matrix techniques. One uses SRRSC for the HSS construction and the other
uses a randomized HSS construction algorithm which is first introduced in~\cite{rand-hss}.
For the later one, we propose a method to estimate the HSS rank
by using the function approximation theory.
The main point is using the rank-structured matrix techniques to update the eigenvectors.
Roughly speaking,  the worst case complexity of ADCs is reduced to $O(N^2r)$ for an $N\times N$
symmetric tridiagonal matrix instead of $O(N^3)$, where $r$ is a modest number which depends on
the property of the tridiagonal matrix.
We implement ADCs in parallel including the HSS construction and HSS matrix multiplication algorithms, and
compare them with the multithreaded Intel MKL library.
For some matrices of large dimensions with few deflations, our ADC algorithms can be more than 6x times faster than the DC algorithm in MKL.

\section*{Acknowledgement}

The authors would like to thank Ming Gu for valuable suggestions and Ziyang Mao,
Lihua Chi, Yihui Yan, Xu Han, Xinbiao Gan
and Qingfeng Hu for some helpful discussions.
The authors also thank the referee for their valuable comments which greatly improve the presentation of this paper.
This work is partial supported by National Natural Science Foundation of China (No. 11401580, 611330005 and 91430218),
and 863 Program of China under grant 2012AA01A301.


\end{document}